\documentclass[12pt,a4paper]{amsart}
\usepackage{geometry}
\usepackage{tabls}
\usepackage{url}
\usepackage[utf8]{inputenc}
\usepackage{amsfonts}
\usepackage{amssymb}
\usepackage{mathrsfs}
\usepackage{amsmath}
\usepackage{amsthm}
\usepackage{amscd}
\usepackage{fancyhdr}
\usepackage{enumerate}
\usepackage{paralist}
\usepackage{graphicx}
\usepackage{cite}

\usepackage{footmisc}
\numberwithin{equation}{section}

\theoremstyle{plain}

\newtheorem{Thm}[equation]{Theorem}
\newtheorem{lem}[equation]{Lemma}
\newtheorem{prop}[equation]{Proposition}

\begin{document}

\title{Remark on a symmetric zeta function}

\author{Jiangtao Li}

\email{lijiangtao@csu.edu.cn}
\address{Jiangtao Li \\ School of Mathematics and Statistics, Central South University, Hunan Province, China}

\begin{abstract}
In this paper we  define a symmetric zeta function. We show that it can be analytically continued to a meromorphic function on $\mathbb{C}^3$ with only simple poles at some special hyperplanes. We also calculate the value of a multiple residue at one special point.  For a divergent multiple series, which can be viewed as the value of the symmetric zeta function at the point $(1,1,1)$, we give a very detailed analysis on its growth and we relate it to the classical Euler constant.
      \end{abstract}

\let\thefootnote\relax\footnotetext{
2020 $\mathnormal{Mathematics} \;\mathnormal{Subject}\;\mathnormal{Classification}$. 11M32.\\
$\mathnormal{Keywords:}$  Riemann zeta function, Multiple zeta values. }

\maketitle

\section{Introduction}\label{int}

The Riemann Zeta Function is defined by  $$\zeta(s)=\sum_{n\geq1}\frac{1}{n^s}.$$
It is convergent when $\mathrm{Re}\;(s)>1$. What is more, 
    $\zeta(s)$ can be analytically continued  a meromorphic function on $\mathbb{C}$ (with  only one simple pole at $s=1$). 
    The multiple zeta function is defined by 
     $$\zeta(s_1, s_2,\cdots s_r)=\sum_{0<n_1<n_2<\cdots<n_r}\frac{1}{n_1^{s_1}n_2^{s_2}\cdots n_r^{s_r}}.$$  
     Zhao \cite{zhao} showed that the multiple zeta function can be analytically continued to  a meromorphic function on $\mathbb{C}^r$ (with  only simple poles at some special hyperplanes). 
                      
      In this paper, we define the following symmetric zeta function:
      \[
      \zeta_{\mathfrak{S}}(s_1,s_2,s_3)=\sum_{n_1,n_2,n_3\geq1}\frac{1}{(n_1+n_2)^{s_1}(n_2+n_3)^{s_2}(n_3+n_1)^{s_3}}.
      \]
     For the symmetric zeta function $\zeta_{\mathfrak{S}}(s_1,s_2,s_3)$, we have:
    \begin{Thm}\label{an}
    The symmetric zeta function $\zeta_{\mathfrak{S}}(s_1,s_2,s_3)$ is convergent when $$\mathrm{Re}(s_1),\mathrm{Re}(s_2),\mathrm{Re}(s_3)>1.$$ It  can be analytically
continued to a meromorphic function on all of $\mathbb{C}^3$, with  simple poles at   the following hyperplanes
$$s_1+s_2=2-k_1, k_1\geq 1, s_2+s_3=2-k_2, k_2\geq 1, s_3+s_1=2-k_3, k_3\geq 1,$$
\[
s_1+s_2+s_3=4-k_4, k_4\geq 1.
\]
  \end{Thm}
  
    By Theorem \ref{an}, now one  can view $\zeta_{\mathfrak{S}}(s_1,s_2,s_3)$  as a meromorphic function on $\mathbb{C}^3$.  Since all its poles are simple, we are interested in the residues of $\zeta_{\mathfrak{S}}(s_1,s_2,s_3)$ at  the simple poles. Actually, we have:
    \begin{Thm}\label{res}
    \[\lim_{(s_1,s_2,s_3)\rightarrow (1,1,1)}(s_1+s_2+s_3-3)\zeta_{\mathfrak{S}}(s_1,s_2,s_3)=\frac{\pi^2}{4},\]
              \end{Thm}

     Cyclotomic multiple zeta values are defined by
      \[
      \zeta\binom{k_1,\cdots,k_r}{\epsilon_1,\cdots,\epsilon_r}=\sum_{0<n_1<\cdots<n_r}\frac{\epsilon_1^{n_1}\cdots \epsilon_r^{n_r}}{n_1^{k_1}\cdots n_r^{k_r}}, (k_r,\epsilon_r)\neq(1,1),
      \]
      where $\epsilon_i^N=1, \forall \,1\leq i\leq r$ for some $N\geq 1$.
      The condition  $(k_r,\epsilon_r)\neq(1,1)$ ensures that the above multiple series are convergent.    Cyclotomic multiple zeta values in cases $\epsilon_1=\cdots=\epsilon_r=1$ are called multiple zeta values.

     By the main results of \cite{ter}, it follows that if 
     \[
     (s_1,s_2,s_3)\in\mathbb{N}^3,\;s_1+s_2+s_3\geq4,\;s_1s_2s_3>0,
     \]     
     then $\zeta_{\mathfrak{S}}(s_1,s_2,s_3)$ is a $\mathbb{Q}[\mu_N]$-linear combination of cyclotomic multiple zeta values for some N, where 
     \[
     \mu_N=\{\epsilon\in\mathbb{C}\;|\;\epsilon^N=1\}.
     \]  
     
     By Theorem \ref{an}, the multiple series
     \[
     \sum_{n_1,n_2,n_3\geq1}\frac{1}{(n_1+n_2)(n_2+n_3)(n_3+n_1)}     
     \]
     is divergent. The following  theorem give a more detailed analysis of the above divergent series.

     \begin{Thm}\label{div}
    For $k\geq 1$, define $f(k)$ as 
     \[
     f(k):=\sum_{1\leq m_i\leq k}\frac{1}{(m_1+m_2)(m_2+m_3)(m_1+m_3)}-\mathop{\int}_{[1,k+1]^3}\frac{dxdydz}{(x+y)(y+z)(z+x)},
     \]
     then
     \[
     \lim_{k\rightarrow +\infty}f(k)=\frac{3}{2}\zeta(2)\gamma+\frac{1}{8}\zeta(3),\]
     where $\gamma$ is the Euler constant which are defined by 
     \[
     \gamma:=\lim_{k\rightarrow +\infty}\left(\sum_{1\leq n\leq k}\frac{1}{n}-\mathrm{log}\;k \right)=\lim_{k\rightarrow +\infty}\left(\sum_{1\leq n\leq k}\frac{1}{n}-\int_1^{k+1}\frac{dt}{t} \right).
     \]
   \end{Thm}
    The limit of $f(k)$ can be viewed as a high dimension analogue of the Euler constant $\gamma$. Theorem \ref{div} gives a detailed analysis of this high dimension analogue.
     
       \section{Analytic continuation of the symmetric zeta function}
       In this section we firstly show that the symmetric zeta function is convergent for $\mathrm{Re}(s_1), \mathrm{Re}(s_2),\mathrm{Re}(s_3)>1$.
       The analytic continuation of $\zeta_{\mathfrak{S}}(s_1,s_2,s_3)$ is similar to the analytic continuation of $\zeta(s)$. We review the theory of  Riemann zeta function shortly in this section. The reference is \cite{zhao}. 
       \subsection{Convergence of the symmetric multiple series}
       For $\mathrm{Re}(s_1), \mathrm{Re}(s_2),\mathrm{Re}(s_3)>1$, denote by
       \[
       s_1=\sigma_1+it_1, s_2=\sigma_2+it_2, s_3=\sigma_3+it_3,
       \]
       where $\sigma_1,\sigma_2,\sigma_3,t_1,t_2,t_3\in\mathbb{R}$. Then we have
       \[
       |\zeta_{\mathfrak{S}}(s_1,s_2,s_3)|\leq \sum_{n_1,n_2,n_3\geq1}\frac{1}{(n_1+n_2)^{\sigma_1}(n_2+n_3)^{\sigma_2}(n_3+n_1)^{\sigma_3}}< \sum_{n_1,n_2,n_3\geq1}\frac{1}{n_1^{\sigma_1}n_2^{\sigma_2}n_3^{\sigma_3}}   .         \]
  So  $\zeta_{\mathfrak{S}}(s_1,s_2,s_3)$ is convergent for $\mathrm{Re}(s_1),\mathrm{Re}(s_2),\mathrm{Re}(s_3)>1$.    
       
       \subsection{Analytic continuation of Riemann zeta function}
       For $\mathrm{Re}\,(s)>0$, the Gamma function $\Gamma(s)$ is defined by 
       \[
       \Gamma(s)=\int^{+\infty}_{0}e^{-t}t^{s-1}dt.
       \]
       It is easy to check that $\Gamma(s+1)=s\Gamma(s)$ for $\mathrm{Re}\,(s)>0$. As a result, 
       \[
       \Gamma(s)=\frac{\Gamma(s+1)}{s}=\cdots=\frac{\Gamma(s+n)}{s(s+1)\cdots (s+n-1)},
       \]
       for $\mathrm{Re}\,(s)>0$. The above formula can be used to continue $\Gamma(s)$ to a meromorphic function over $\mathbb{C}$, with simple poles at $s=0,-1,-2,\cdots,-n,\cdots$. One can check that $$\mathrm{Res}\,\big{|}_{s=0} \Gamma(s)=1,\,\mathrm{Res}\,\big{|}_{s=-n} \Gamma(s)=\frac{(-1)^{n}}{n!},n\geq 1.$$
       
       For $\mathrm{Re}(s)>1,n\geq 1$, it is easy to check that 
        \[
       \frac{\Gamma(s)}{n^s}=\int^{+\infty}_{0}e^{-nt}t^{s-1}dt.
       \]       
       Thus 
       \[
       \Gamma(s)\zeta(s)=\sum_{n\geq1}\int^{+\infty}_{0}e^{-nt}t^{s-1}dt=\int^{+\infty}_0\frac{t^{s-1}dt}{e^t-1}.
              \]
              Now we have 
              \[
              \Gamma(s)\zeta(s)=\int^1_0t^{s-2}\frac{tdt}{e^t-1}+\int^{+\infty}_1\frac{t^{s-1}dt}{e^t-1}.
              \]
       Denote by $$M(s)=\int^1_0t^{s-2}\frac{tdt}{e^t-1}, \,R(s)=\int^{+\infty}_1\frac{t^{s-1}dt}{e^t-1},$$   
       it is easy to see that $R(s)$ is a holomorphic functions over $\mathbb{C}$. For $M(s)$, we need the following lemma:
       \begin{lem}\label{tay}
    For a fixed $\epsilon>0$, if $\varphi(t)$ is a smooth (infinitely differentiable)  bounded function on $(-\epsilon,1)$, then the function 
    \[
    M_{\varphi}(s)=\int^1_0t^{s-2}\varphi(t)dt
    \]
    can be analytically continued to a meromorphic function on $\mathbb{C}$, with simple poles at $s=2-k,k\geq 1$, and $\mathrm{Res}_{s=2-k}M_{\varphi}(s)=\frac{\varphi^{(k-1)}(0)}{(k-1)!}$.
    \end{lem}
     \noindent{\bf Proof:} For $\mathrm{Re}(s)>1$ and any fixed $k$, one can that
     \[
     M_{\varphi}(s)=\frac{\varphi(0)}{s-1}+\sum_{n=1}^k\frac{\varphi^{(n)}(0)}{n!}\frac{1}{s+n-1}+\int^1_0t^{s-2}\left[ \varphi(t)-\varphi(0)-\sum_{n=1}^k\frac{\varphi^{(n)}(0)}{n!}t^n  \right]dt.
     \]     
     Since $\varphi(t)$ is an infinitely differentiable function on $(-\epsilon,1)$,  we have 
     \[
     R_{\varphi}(t)=\varphi(t)-\varphi(0)-\sum_{n=1}^k\frac{\varphi^{(n)}(0)}{n!}t^n=\mathrm{O}(t^{k+1}).
     \]
     Thus the function 
     \[
     \int^1_0t^{s-2}R_{\varphi}(t)dt
     \]
     can be analytically continued to a holomorphic function on  $\{s\in\mathbb{C}\,|\,  \mathrm{Re}(s)>-k\}$.
    As a result, the lemma is proved.
$\hfill\Box$\\

Letting $\varphi(t)=\frac{t}{e^t-1}$, by Lemma \ref{tay} and comparing the poles of $M_{\varphi}(s)$ and $\Gamma(s)$, one can deduce that $\zeta(s)$ can be analytically continued to a meromorphic function on $\mathbb{C}$ with an only simple pole at $s=1$. Moreover, $\mathrm{Res}_{s=1}\zeta(s)=1$.

\subsection{Analytic continuation of the symmetric zeta function}\label{anly}
For $\mathrm{Re}(s_i)>1,n_i\geq 1$, we have 
\[
\frac{\Gamma(s_1)}{(n_1+n_2)^{s_1}}=\int^{+\infty}_0e^{-(n_1+n_2)t_1}t_1^{s_1-1}dt_1,\]
\[\frac{\Gamma(s_2)}{(n_2+n_3)^{s_2}}=\int^{+\infty}_0e^{-(n_2+n_3)t_2}t_2^{s_2-1}dt_2,\;\;\frac{\Gamma(s_3)}{(n_3+n_1)^{s_3}}=\int^{+\infty}_0e^{-(n_3+n_1)t_3}t_3^{s_3-1}dt_3.\]
Thus for $\mathrm{Re}(s_i)>1$, 
it follows that         
\[
\begin{split}
&\;\;\;\;\,\Gamma(s_1)\Gamma(s_2)\Gamma(s_3)\zeta_{\mathfrak{S}}(s_1,s_2,s_3)\\
&=\sum_{n_1,n_2,n_3\geq 1}\mathop{\iiint}_{[0,+\infty)^3} e^{-(n_1+n_2)t_1}e^{-(n_2+n_3)t_2}e^{-(n_3+n_1)t_3}  t_1^{s_1-1}t_2^{s_2-1}t_3^{s_3-1} dt_1dt_2dt_3   \\
&=\sum_{n_1,n_2,n_3\geq 1} \mathop{\iiint}_{[0,+\infty)^3} e^{-n_1(t_3+t_1)}e^{-n_2(t_1+t_2)}e^{-n_3(t_2+t_3)}     t_1^{s_1-1}t_2^{s_2-1}t_3^{s_3-1} dt_1dt_2dt_3   \\
&= \mathop{\iiint}_{[0,+\infty)^3}\frac{  t_1^{s_1-1}t_2^{s_2-1}t_3^{s_3-1} dt_1dt_2dt_3 }{ (e^{t_3+t_1}-1)(e^{t_1+t_2}-1)(e^{t_2+t_3}-1)}  \\
&= \mathop{\iiint}_{[0,+\infty)^3}\frac{t_1^{s_1-1}t_2^{s_2-1}t_3^{s_3-1}}{(t_3+t_1)(t_1+t_2)(t_2+t_3)}\frac{(t_3+t_1)(t_1+t_2)(t_2+t_3)}{(e^{t_3+t_1}-1)(e^{t_1+t_2}-1)(e^{t_2+t_3}-1)}dt_1dt_2dt_3.
\end{split}
\]        
For  convenience, we let $\Theta(s_1,s_2,s_3)=\Gamma(s_1)\Gamma(s_2)\Gamma(s_3)\zeta_{\mathfrak{S}}(s_1,s_2,s_3)$ and 
\[
\varphi(t_1,t_2,t_3)=\frac{(t_3+t_1)(t_1+t_2)(t_2+t_3)}{(e^{t_3+t_1}-1)(e^{t_1+t_2}-1)(e^{t_2+t_3}-1)}.
\]
We have 
\[
\begin{split}
&\;\;\;\;\,\Theta(t_1,t_2,t_3)\\
&=\mathop{\iiint}_{[0,+\infty)^3}\frac{t_1^{s_1-1}t_2^{s_2-1}t_3^{s_3-1}}{(t_3+t_1)(t_1+t_2)(t_2+t_3)}\varphi(t_1,t_2,t_3)dt_1dt_2dt_3 \\
&=M_1(s_1,s_2,s_3)+M_2(s_1,s_2,s_3)+R_1(s_1,s_2,s_3)+R_2(s_1,s_2,s_3),   \\
\end{split}
\]
where 
\[
M_1(s_1,s_2,s_3)=\mathop{\iiint}_{[0,1]^3}\frac{t_1^{s_1-1}t_2^{s_2-1}t_3^{s_3-1}}{(t_3+t_1)(t_1+t_2)(t_2+t_3)}\varphi(t_1,t_2,t_3)dt_1dt_2dt_3,
\]
\[
\begin{split}
&\;\;\;\;\,M_2(s_1,s_2,s_3)\\
&=\left(\mathop{\iiint}_{[0,1]^2\times[1,+\infty)}+\mathop{\iiint}_{[0,1]\times[1,+\infty)\times [0,1]}+\mathop{\iiint}_{[1,+\infty)\times[0,1]^2}\right)\frac{t_1^{s_1-1}t_2^{s_2-1}t_3^{s_3-1}\varphi(t_1,t_2,t_3)dt_1dt_2dt_3}{(t_3+t_1)(t_1+t_2)(t_2+t_3)},\\
\end{split}
\]
\[
\begin{split}
&\;\;\;\;\,R_1(s_1,s_2,s_3)\\
&=\left(\mathop{\iiint}_{[0,1]\times[1,+\infty)^2}+\mathop{\iiint}_{[1,\infty)\times[0,1]\times [1,\infty)}+\mathop{\iiint}_{[1,+\infty)^2\times[0,1]}\right)\frac{t_1^{s_1-1}t_2^{s_2-1}t_3^{s_3-1}\varphi(t_1,t_2,t_3)dt_1dt_2dt_3}{(t_3+t_1)(t_1+t_2)(t_2+t_3)},\\
\end{split}
\]
\[
R_2(s_1,s_2,s_3)=\mathop{\iiint}_{[1,+\infty)^3}\frac{t_1^{s_1-1}t_2^{s_2-1}t_3^{s_3-1}}{(t_3+t_1)(t_1+t_2)(t_2+t_3)}\varphi(t_1,t_2,t_3)dt_1dt_2dt_3.
\]

It is clear that $R_2(s_1,s_2,s_3)$ is a holomorphic function on $\mathbb{C}^3$. By Lemma \ref{tay}, it follows that $R_1(s_1,s_2,s_3)$ can be analytically continued to a meromorphic function on $\mathbb{C}^3$, with simple poles at the following co-dimension  one hyperplanes $$s_1=1-k_1,k_1\geq 1, s_2=1-k_2, k_2\geq 1,s_3=1-k_3, k_3\geq 1.$$

The analysis of $M_2(s_1,s_2,s_3)$ is a little bit tricky. For convenience, define 
\[
f(s_1,s_2,s_3)=\mathop{\iiint}_{[0,1]^2\times[1,+\infty)}\frac{t_1^{s_1-1}t_2^{s_2-1}t_3^{s_3-1}}{(t_3+t_1)(t_1+t_2)(t_2+t_3)}\varphi(t_1,t_2,t_3)dt_1dt_2dt_3.
\]
By the symmetry of $\frac{\varphi(t_1,t_2,t_3)}{(t_3+t_1)(t_1+t_2)(t_2+t_3)}$, it is easy to check that 
\[
M_2(s_1,s_2,s_3)=f(s_1,s_2,s_3)+f(s_3,s_1,s_2)+f(s_2,s_3,s_1).
\]
Now we focus on $f(s_1,s_2,s_3)$. It can be split up into two parts, i.e.
\[
f(s_1,s_2,s_3)=\left(\mathop{\iiint}_{\Delta_1\times[1,+\infty)}+\mathop{\iiint}_{\Delta_2\times[1,+\infty)}  \right)\frac{t_1^{s_1-1}t_2^{s_2-1}t_3^{s_3-1}}{(t_3+t_1)(t_1+t_2)(t_2+t_3)}\varphi(t_1,t_2,t_3)dt_1dt_2dt_3,
\]
where 
\[
\Delta_1=\big{\{}(t_1,t_2)\in [0,1]^2\,|\,0\leq t_1\leq t_2\leq 1\big{\}},\Delta_2=\big{\{}(t_1,t_2)\in [0,1]^2\,|\,0\leq t_2\leq t_1\leq 1\big{\}}.
\]

On $\Delta_1$, we use the changing of variables: $t_1=uv,t_2=v, (u,v)\in [0,1]^2$. On $\Delta_2$, we use the changing of variables: $t_1=u,t_2=uv, (u,v)\in [0,1]^2$.  Then we have 
\[
\begin{split}
&\;\;\;\;\;f(s_1,s_2,s_3)\\
&=\int^{+\infty}_1dt_3\mathop{\iint}_{[0,1]^2}\frac{u^{s_1-1}v^{s_1+s_2-2}t_3^{s_3-1}}{(t_3+uv)(1+u)(v+t_3)}\varphi(uv,v,t_3)dudv\\
&+\int^{+\infty}_1dt_3\mathop{\iint}_{[0,1]^2}\frac{u^{s_1+s_2-2}v^{s_2-1}t_3^{s_3-1}}{(t_3+u)(1+v)(uv+t_3)}\varphi(u,uv,t_3)dudv.\\
\end{split}
\]

By Lemma \ref{tay},  $f(s_1,s_2,s_3)$ can be analytically continued to a meromorphic function on $\mathbb{C}^3$, with simple poles at the following co-dimension  one hyperplanes $$s_1=1-k_1,k_1\geq 1, s_2=1-k_2, k_2\geq 1, s_1+s_2=2-k_3, k_3\geq 1.$$
Thus $M_2(s_1,s_2,s_3)$ also can be analytically continued to a meromorphic function on $\mathbb{C}^3$, with simple poles at the following co-dimension  one hyperplanes $$s_1=1-k_1,k_1\geq 1, s_2=1-k_2, k_2\geq 1, s_3=1-k_3, k_3\geq 1,$$
$$s_1+s_2=2-k_4, k_4\geq 1, s_3+s_1=2-k_5,k_5\geq 1, s_2+s_3=2-k_6, k_6\geq 1.$$

Denote by
\[
g(s_1,s_2,s_3)=\mathop{\iiint}_{0<t_1<t_2<t_3<1}\frac{t_1^{s_1-1}t_2^{s_2-1}t_3^{s_3-1}}{(t_3+t_1)(t_1+t_2)(t_2+t_3)}\varphi(t_1,t_2,t_3)dt_1dt_2dt_3,
\]
then
\[
\begin{split}
&\;\;\;\;M_1(s_1,s_2,s_3)\\
&=g(s_1,s_2,s_3)+g(s_2,s_3,s_1)+g(s_3,s_1,s_2)+g(s_2,s_1,s_3)+g(s_1,s_3,s_2)+g(s_3,s_2,s_1).   \\
\end{split}
\]
By changing of variables $t_1=uvw,t_2=vw,t_3=w$, $g(s_1,s_2,s_3)$ is reduced to
\[
g(s_1,s_2,s_3)=\mathop{\iiint}_{[0,1]^3}\frac{u^{s_1-1}v^{s_1+s_2-2}w^{s_1+s_2+s_3-4}}{(1+uv)(1+u)(1+v)}\varphi(uvw,vw,w)dudvdw. \tag{1}
\]

By Lemma \ref{tay},  $g(s_1,s_2,s_3)$ can be analytically continued to a meromorphic function on $\mathbb{C}^3$, with simple poles at the following co-dimension  one hyperplanes $$s_1=1-k_1,k_1\geq 1, s_1+s_2=2-k_2, k_2\geq 1, s_1+s_2+s_3=4-k_3, k_3\geq 1.$$
So $M_1(s_1,s_2,s_3)$ can be analytically continued to a meromorphic function on $\mathbb{C}^3$, with simple poles at the following co-dimension  one hyperplanes $$s_1=1-k_1,k_1\geq 1, s_2=1-k_2,k_2\geq 1, s_3=1-k_3, k_3\geq 1,    $$
$$s_1+s_2=2-k_4, k_4\geq 1, s_2+s_3=2-k_5, k_5\geq 1, s_3+s_1=2-k_6, k_6\geq 1,$$
$$ s_1+s_2+s_3=4-k_7, k_7\geq 1.$$

From the above analysis, we find that the symmetric zeta function 
$$\zeta_{\mathfrak{S}}(s_1,s_2,s_3)$$ can be analytically continued to a meromorphic function on $\mathbb{C}^3$, with simple poles at the following co-dimension  one hyperplanes $$s_1+s_2=2-k_1, k_1\geq 1, s_2+s_3=2-k_2, k_2\geq 1, s_3+s_1=2-k_3, k_3\geq 1,$$
$$s_1+s_2+s_3=4-k_4, k_4\geq 1.$$
As a result, Theorem \ref{an} is proved.

              \section{The residue at the simple pole}
              In this section we compute the value of the multiple residue of the symmetric zeta function at the point $(1,1,1)$.

              From the analysis in Section \ref{anly}, it follows that 
              \[
              \begin{split}
             &\;\;\;\; \lim_{(s_1,s_2,s_3)\rightarrow(1,1,1)}(s_1+s_2+s_3-3)\zeta_{\mathfrak{S}}(s_1,s_2,s_3)\\
             &=\lim_{(s_1,s_2,s_3)\rightarrow(1,1,1)} (s_1+s_2+s_3-3)\frac{M_1(s_1,s_2,s_3)}{\Gamma(s_1)\Gamma(s_2)\Gamma(s_3)}\\
             &=\lim_{(s_1,s_2,s_3)\rightarrow(1,1,1)} (s_1+s_2+s_3-3)M_1(s_1,s_2,s_3)\\
             &=\lim_{(s_1,s_2,s_3)\rightarrow(1,1,1)}6(s_1+s_2+s_3-3)g(s_1,s_2,s_3).\\
                        \end{split}
              \]
              By Lemma \ref{tay} and formula $(1)$, we have
              \[
              \lim_{(s_1,s_2,s_3)\rightarrow(1,1,1)}(s_1+s_2+s_3-3)g(s_1,s_2,s_3)=6\mathop{\iint}_{[0,1]^2}\frac{dudv}{(1+u)(1+v)(1+uv)}.
                            \]
                            
                            In a word, Theorem \ref{res} is reduced to the following statement.
       \begin{Thm}\label{inte}
      \[ \mathop{\iint}_{[0,1]^2}\frac{dudv}{(1+u)(1+v)(1+uv)}=\frac{\pi^2}{24}.\]
              \end{Thm}
       
                \noindent{\bf Proof:}  It is easy to check that
                \[
                \frac{1}{(1+u)(1+v)(1+uv)}=\frac{1}{v(1+v)}\cdot \frac{1}{(1+u)(\frac{1}{v}+u)}=\frac{1}{1-v^2}\left(\frac{1}{1+u}-\frac{1}{\frac{1}{v}+u}    \right).
                \]
                Thus we have 
                \[
                \begin{split}
                &\;\;\;\;\mathop{\iint}_{[0,1]^2}\frac{dudv}{(1+u)(1+v)(1+uv)}      \\
                &=\mathop{\mathrm{lim}}_{\epsilon\rightarrow 0^+}\mathop{\int}_{0}^{1-\epsilon}dv\mathop{\int}^1_0 \frac{1}{1-v^2}\left(\frac{1}{1+u}-\frac{1}{\frac{1}{v}+u}    \right)du           \\
                &=\mathop{\mathrm{lim}}_{\epsilon\rightarrow 0^+}\mathop{\int}_{0}^{1-\epsilon} \frac{1}{1-v^2} \left[ \mathrm{log}\;2-\mathrm{log}\;(1+v)        \right]       dv   \\
                &=\frac{1}{2}\cdot \mathop{\mathrm{lim}}_{\epsilon\rightarrow 0^+}\mathop{\int}_{0}^{1-\epsilon}\left( \frac{1}{1+v} +\frac{1}{1-v}\right)\left[ \mathrm{log}\;2-\mathrm{log}\;(1+v)        \right]       dv   \\  
                &=\frac{1}{2}(\mathrm{I}+\mathrm{II}).   \end{split}
     \]
     Here 
     \[
     \mathrm{I}=\mathop{\mathrm{lim}}_{\epsilon\rightarrow 0^+}\mathop{\int}_{0}^{1-\epsilon} \frac{1}{1+v} \left[ \mathrm{log}\;2-\mathrm{log}\;(1+v)        \right]       dv       \]
     and 
       \[
     \mathrm{II}=\mathop{\mathrm{lim}}_{\epsilon\rightarrow 0^+}\mathop{\int}_{0}^{1-\epsilon} \frac{1}{1-v} \left[ \mathrm{log}\;2-\mathrm{log}\;(1+v)        \right]       dv .      \]        
     
     For $\mathrm{I}$, one has
     \[
     \begin{split}
     &\;\;\;\;\mathrm{I}           \\
     &=  \mathop{\int}_{0}^{1} \frac{1}{1+v} \left[ \mathrm{log}\;2-\mathrm{log}\;(1+v)        \right]       dv        \\
     &=\left[ \mathrm{log}\,2\mathrm{log}(1+v)-\frac{1}{2}\mathrm{log}^2(1+v)\right]\Bigg{|}^{v=1}_{v=0}\\
     &=\frac{1}{2}\mathrm{log}^2 2.
     \end{split}
     \]       
     For $\mathrm{II}$, one has
     \[
     \begin{split}
     &\;\;\;\;\mathrm{II}\\
     &=  \mathop{\mathrm{lim}}_{\epsilon\rightarrow 0^+}\mathop{\int}_{\epsilon}^{1} \frac{1}{v} \left[ \mathrm{log}\;2-\mathrm{log}\;(2-v)        \right]       dv            
     = \mathop{\mathrm{lim}}_{\epsilon\rightarrow 0^+}\mathop{\int}_{\epsilon}^{1}- \frac{1}{v}\mathrm{log}\; (1-\frac{v}{2})dv \\
     &= \mathop{\mathrm{lim}}_{\epsilon\rightarrow 0^+}\mathop{\int}_{\epsilon}^{1} \sum_{l\geq 1}\frac{v^{l-1}}{2^l\cdot l}  dv
     =\sum_{l\geq 1}\frac{1}{2^l}\frac{1}{l^2}.
     \end{split}
     \]
     By the formula (7.7) in \cite{bbb}, one has 
     \[
    \mathrm{II}= \mathrm{Li}_2\left(\frac{1}{2}\right)=\sum_{l\geq 1}\frac{1}{2^l}\cdot \frac{1}{l^2}=\frac{\pi^2}{12}-\frac{1}{2}\mathrm{log}^2 2.
     \]
     As a result, we have 
    \[ \mathop{\iint}_{[0,1]^2}\frac{dudv}{(1+u)(1+v)(1+uv)}=\frac{1}{2}(\mathrm{I}+\mathrm{II})=\frac{\pi^2}{24}.
    \]                                       $\hfill\Box$\\

       \section{Analysis of the divergent series}
                     In this section we deal with the divergent multiple series
                     \[
                     \sum_{m_1,m_2,m_3\geq 1}\frac{1}{(m_1+m_2)(m_2+m_3)(m_3+m_1)}.
                     \]
                     Firstly we show that a modified version of the above series is convergent to a real number. 
                     Then we will discuss a partial sum of finite terms of the divergent multiple series.
                     Lastly we will establish the relation between this  real number and the famous Euler constant.
                     
    \begin{prop}\label{1}
    For $k\geq 1$, define $f(k)$ as 
     \[
     f(k):=\sum_{1\leq m_i\leq k}\frac{1}{(m_1+m_2)(m_2+m_3)(m_1+m_3)}-\mathop{\iiint}_{[1,k+1]^3}\frac{dxdydz}{(x+y)(y+z)(z+x)},
     \]
     then the limit $\lim\limits_{k\rightarrow +\infty}f(k)$
     exists.
    \end{prop}                 
  \noindent{\bf Proof:} One can rewrite $f(k)$ as
  \[
  \begin{split}
 & f(k)=\sum_{1\leq m_i\leq k}\mathop{\iiint}_{[0,1]^3}
  \Bigg{(}\frac{1}{(m_1+m_2)(m_2+m_3)(m_3+m_1)} \\
  &\;\;\;\;\;\;\;\;\;\;\;\;\;\; -\frac{1}{(m_1+m_2+x+y)(m_2+m_3+y+z)(m_3+m_1+z+x)}    \Bigg{)}dxdydz.\\
  \end{split}
  \]
 From the above expression, it is clear that $$f(k)<f(k+1),\forall\;k\geq 1.$$
  Thus it suffices to prove that $f(k)$ is a bounded sequence.
  
  For any fixed $(x,y,z)\in[0,1]^3$, define $F_{m_1,m_2,m_3}(t)$, $t\in {[0,1]}$ as 
  \[
  F_{m_1,m_2,m_3}(t)=\frac{1}{[m_1+m_2+t(x+y)][m_2+m_3+t(y+z)][m_3+m_1+t(z+x)]}.
  \]
  Then  
  \[
  \begin{split}
  & \frac{1}{(m_1+m_2)(m_2+m_3)(m_3+m_1)} \\
  &\;\;\;\;\;\;\;\;\;\;\;\;\;\;\;\;\;\;\;\;\;\;\;\;\;\; -\frac{1}{(m_1+m_2+x+y)(m_2+m_3+y+z)(m_3+m_1+z+x)} \\
  & =F(0)-F(1)=-F^\prime(t_{x,y,z})     \\
  \end{split}
  \]
  for some $t_{x,y,z}\in (0,1)$. Since 
  \[
  \begin{split}
  &\;\;\;\;F^\prime(t)\\
  &=\frac{-1}{[m_1+m_2+t(x+y)][m_2+m_3+t(y+z)][m_3+m_1+t(z+x)]}       \\
  &\left[\frac{x+y}{m_1+m_2+t(x+y)} + \frac{y+z}{m_2+m_3+t(y+z)} +\frac{z+x}{m_3+m_1+t(z+x)} \right],                                 \\
  \end{split}
    \]
    we have 
    \[
    \begin{split}
     & 0< \mathop{\iiint}_{[0,1]^3} \Bigg{(}\frac{1}{(m_1+m_2)(m_2+m_3)(m_3+m_1)} \\
  &\;\;\;\;\;\;\;\;\;\;\;\;\;\; -\frac{1}{(m_1+m_2+x+y)(m_2+m_3+y+z)(m_3+m_1+z+x)}    \Bigg{)}dxdydz\\
  \end{split}
  \]
  \[
  \begin{split}
  &<\frac{2}{(m_1+m_2)^2(m_2+m_3)(m_3+m_1)}+\frac{2}{(m_1+m_2)(m_2+m_3)^2(m_3+m_1)} +\\
  &\;\;\;\;\;\frac{2}{(m_1+m_2)(m_2+m_3)(m_3+m_1)^2} 
  \end{split}    \]
  By Section \ref{anly},  $\zeta_{\mathfrak{S}}(2,1,1)=\zeta_{\mathfrak{S}}(1,2,1)=\zeta_{\mathfrak{S}}(1,1,2)$ are convergent. So $$0<f(k)<6\zeta_{\mathfrak{S}}(2,1,1), \forall\; k\geq 1.$$      $\hfill\Box$\\   
  
  Now we treat the triple sum and the triple integral separately.  
  
  \subsection{Triple harmonic sum}
  
  The following  lemmas will play  important roles in our analysis.
  
   \begin{lem}\label{pow}
  For $l,n\geq 1$, there is an $M$ (which is independent of $l$ and $n$), such that 
  \[
 \; \Bigg{|}\sum_{1\leq m\leq n}m^l-\frac{n^{l+1}}{l+1}\Bigg{|}\leq Mn^l, \Bigg{|}\sum_{1\leq m< n}m^l-\frac{n^{l+1}}{l+1}\Bigg{|}\leq Mn^l  \]
and 
 \[
  \Bigg{|}\sum_{1\leq m< n}m^l-\frac{n^{l+1}}{l+1}+\frac{n^l}{2}\Bigg{|}= \Bigg{|}\sum_{1\leq m\leq  n}m^l-\frac{n^{l+1}}{l+1}-\frac{n^l}{2}\Bigg{|}\leq Mln^{l-1}.\]
  \end{lem}
   \noindent{\bf Proof:}  By Theorem $5.1$ (the Euler–Maclaurin summation formula) in Chapter $5$, \cite{kan}, we have
   \[
   \sum_{1\leq m\leq n}m^l= \sum_{0\leq m\leq n}m^l =\frac{n^{l+1}}{l+1}+\frac{n^l}{2}+l\int^n_0B_1(x-[x])x^{l-1}dx
     \]
     and 
      \[
   \sum_{1\leq m\leq n}m^l= \sum_{0\leq m\leq n}m^l =\frac{n^{l+1}}{l+1}+\frac{n^l}{2}+\frac{l}{12}n^{l-1}-\frac{l(l-1)}{2}\int^n_0B_2(x-[x])x^{l-2}dx.
     \]
     
     Here $[x]$ means the integer part of $x$ and $B_1(x), B_2(x)$ are two of   the Bernoulli polynomials. As we don't use the definition of Bernoulli polynomials in our paper, we will not give their detailed definition here. One can consult \cite{kan} for more information.
     
     Denote by  $C_1=\mathop{\mathrm{sup}}\limits_{x\in \mathbb{R}} |B_1(x-[x])| $, then
     \[
      \Bigg{|}\sum_{1\leq m\leq n}m^l -\frac{n^{l+1}}{l+1}-\frac{n^l}{2}\Bigg{|}\leq C_1l\int^n_0x^{l-1}dx=C_1n^l.
           \]
           Define $M^\prime:=C_1+\frac{1}{2}$, then we have 
           \[
            \Bigg{|}\sum_{1\leq m\leq n}m^l -\frac{n^{l+1}}{l+1}\Bigg{|}\leq M^\prime n^l.           \]
            Denote by $M_1=M^\prime+1$, it is clear that 
            \[
             \Bigg{|}\sum_{1\leq m\leq n}m^l-\frac{n^{l+1}}{l+1}\Bigg{|}\leq M_1n^l, \Bigg{|}\sum_{1\leq m< n}m^l-\frac{n^{l+1}}{l+1}\Bigg{|}\leq M_1n^l  .           \]
             
             Similarly, denote by  $C_2=\mathop{\mathrm{sup}}\limits_{x\in \mathbb{R}} |B_2(x-[x])| $, then (assume $l\geq 2$ here)
              \[
      \Bigg{|}\sum_{1\leq m\leq n}m^l -\frac{n^{l+1}}{l+1}-\frac{n^l}{2}-\frac{l}{12}n^{l-1}\Bigg{|}\leq C_2\frac{l(l-1)}{2}\int^n_0x^{l-2}dx<\frac{C_2l}{2}n^{l-1}.
           \]
           Denote by $M_2:=\frac{1}{12}+\frac{C_2}{2}$, we have
           \[
             \Bigg{|}\sum_{1\leq m\leq n}m^l -\frac{n^{l+1}}{l+1}-\frac{n^l}{2}\Bigg{|}\leq M_2ln^{l-1}.
           \]
           
           Beware that the above inequation also holds for $l=1$.  In conclusion, by letting $$M:=max\{M_1,M_2\},$$ the lemma is proved. $\hfill\Box$\\     
           
           \begin{lem}\label{dpow} Let $M$ be the constant in Lemma \ref{pow}. For $l\geq 1,k> 1$, we have 
           \[
          \Bigg{|} \sum_{1\leq m<k}\frac{k^{l+1}-m^{l+1}}{k-m}-\left(1+\frac{1}{2}+\cdots+\frac{1}{l+1}\right)k^{l+1}\Bigg{|}\leq {(M+1)}lk^l.
           \]
           
           \end{lem}
            \noindent{\bf Proof:} By Lemma \ref{pow}, one has
            \[
            \begin{split}
              &\;\;\;\; \Bigg{|} \sum_{1\leq m<k}\frac{k^{l+1}-m^{l+1}}{k-m}-\left(1+\frac{1}{2}+\cdots+\frac{1}{l+1}\right)k^{l+1}\Bigg{|}\\
              &=  \Bigg{|}  \sum_{1\leq m<k}\left(  \sum_{0\leq j\leq l}k^{l-j}m^{j}\right) - \left(1+\frac{1}{2}+\cdots+\frac{1}{l+1}\right)k^{l+1}          \Bigg{|}              \\
               &=  \Bigg{|}  \sum_{ 0\leq j\leq l    }\left(  \sum_{ 1\leq m<k}k^{l-j}m^{j}\right) - \left(1+\frac{1}{2}+\cdots+\frac{1}{l+1}\right)k^{l+1}          \Bigg{|}              \\
                &=  \Bigg{|}  \sum_{ 0\leq j\leq l    }\left(  \sum_{ 1\leq m<k}k^{l-j}m^{j}-\frac{1}{j+1}k^{l+1}\right)   \Bigg{|}  \leq k^l+   \Bigg{|}  \sum_{ 1\leq j\leq l    }\left(  \sum_{ 1\leq m<k}k^{l-j}m^{j}-\frac{1}{j+1}k^{l+1}\right)   \Bigg{|}          \\
                &\leq       k^l+\sum_{ 1\leq j\leq l } \Bigg{|}  \sum_{ 1\leq m<k}k^{l-j}m^{j}-\frac{1}{j+1}k^{l+1} \Bigg{|}  \leq k^l+ lMk^l \leq (M+1)lk^l.  \\
              \end{split}
                          \]
   $\hfill\Box$\\

   \begin{lem}\label{bipow}  For $k>1, p,q\geq 1$, we have 
           \[
          \Bigg{|} \sum_{1\leq m<k}\frac{(k^{p}-m^{p})(k^q-m^q)}{k-m}-\left(\sum_{1\leq n_1\leq p}\frac{1}{n_1}+ \sum_{1\leq n_2\leq q}\frac{1}{n_2}-\sum_{1\leq n_3\leq p+q}\frac{1}{n_3}\right)k^{p+q}  \Bigg{|} \leq 2Mp k^{p+q-1}.
           \]
           \end{lem}
            \noindent{\bf Proof:} By Lemma \ref{pow}, one has
             \[
             \begin{split}
          &\;\;\;\;\Bigg{|} \sum_{1\leq m<k}\frac{(k^{p}-m^{p})(k^q-m^q)}{k-m}-\left(\sum_{1\leq n_1\leq p}\frac{1}{n_1}+ \sum_{1\leq n_2\leq q}\frac{1}{n_2}-\sum_{1\leq n_3\leq p+q}\frac{1}{n_3}\right)k^{p+q}  \Bigg{|} \\
          &=\Bigg{|}\sum_{1\leq m<k}\sum_{1\leq n_1\leq p}k^{p-n_1}m^{n_1-1}(k^q-m^q)-\sum_{1\leq n_1\leq p} \left(\frac{1}{n_1}-\frac{1}{n_1+q} \right)   k^{p+q} \Bigg{|}  \\
          &=  \Bigg{|}\sum_{1\leq n_1\leq p} k^{p-n_1} \left[ \sum_{1\leq m<k}\left(m^{n_1-1}k^q-m^{n_1+q-1}\right)-  \left(\frac{1}{n_1}-\frac{1}{n_1+q} \right)   k^{n_1+q}          \right] \Bigg{|}    \\
          &\leq \sum_{1\leq n_1\leq p}k^{p-n_1}\cdot 2M k^{n_1+q-1}=2Mpk^{p+q-1}.\\
          \end{split}
           \]
            $\hfill\Box$\\

   \begin{lem}\label{ddpow} For $l\geq 1$, $k>1$, we have
   \[
   \begin{split}
  &\Bigg{|} \sum_{1\leq m_1\leq m_2<k}\frac{k^{l+1}-m_1^{l+1}}{(k-m_1)(k-m_2)}\\
  &-\left[ \sum_{0\leq p\leq l}\frac{1}{p+1}\sum_{1\leq m_2<k}\frac{1}{m_2}-\sum_{0\leq p\leq l}\frac{1}{p+1}\left( 1+\frac{1}{2}+\cdots+\frac{1}{p+1}\right)     \right]   k^{l+1} \Bigg{|}\leq \\
  &\leq Mlk^l\sum_{1\leq m<k}\frac{1}{m}+(M+1)(l+1)k^l
    \end{split}
   \]
   \end{lem}  
     \noindent{\bf Proof:}   
     Since 
     \[
     \begin{split}
    &\;\;\;\; \sum_{1\leq m_1\leq m_2<k}\frac{k^{l+1}-m_1^{l+1}}{(k-m_1)(k-m_2)}\\
    &=\sum_{1\leq m_1\leq m_2<k}\sum_{0\leq p\leq l}\frac{k^{l-p}m_1^p}{k-m_2}=\sum_{0\leq p\leq l}k^{l-p}\sum_{1\leq m_1\leq m_2<k}\frac{m_1^p}{k-m_2},\\
    \end{split}
          \]
          by Lemma \ref{pow}, we have 
          \[
          \begin{split}
           &\Bigg{|}\sum_{1\leq m_1\leq m_2<k}\frac{k^{l+1}-m_1^{l+1}}{(k-m_1)(k-m_2)}-\sum_{0\leq p\leq l} \frac{k^{l-p}}{p+1}\sum_{1\leq m_2<k}\frac{m_2^{p+1}}{k-m_2}\Bigg{|}\\
           &\leq \sum_{1\leq p\leq l}k^{l-p}\sum_{1\leq m_2<k}\frac{Mm_2^p}{k-m_2}<Mlk^l\sum_{1\leq m_2<k}\frac{1}{k-m_2}=Mlk^l\sum_{1\leq m<k}\frac{1}{m}.
           \\
                     \end{split}
          \]
          By Lemma \ref{dpow}, we have
          \[
          \begin{split}
          &\;\;\;\;\Bigg{|} \sum_{0\leq p\leq l} \frac{k^{l-p}}{p+1}\sum_{1\leq m_2<k}\frac{m_2^{p+1}}{k-m_2}  \\
           &\;\;\;\;-\left[ \sum_{0\leq p\leq l}\frac{1}{p+1}\sum_{1\leq m_2<k}\frac{1}{m_2}-\sum_{0\leq p\leq l}\frac{1}{p+1}\left( 1+\frac{1}{2}+\cdots+\frac{1}{p+1}\right)     \right]   k^{l+1}\Bigg{|}\\
                     &=\Bigg{|} \sum_{0\leq p\leq l} \frac{k^{l-p}}{p+1}\sum_{1\leq m_2<k}\frac{m_2^{p+1}}{k-m_2}  \\
           &\;\;\;\;- \sum_{0\leq p\leq l}\frac{k^{l-p}}{p+1}\sum_{1\leq m_2<k}\frac{k^{p+1}}{k-m_2}+\sum_{0\leq p\leq l}\frac{1}{p+1}\left( 1+\frac{1}{2}+\cdots+\frac{1}{p+1}\right)       k^{l+1}\Bigg{|}\\           
           &=\Bigg{|}\sum_{0\leq p\leq l}\frac{k^{l-p}}{p+1}\sum_{1\leq m_2<k}\frac{k^{p+1}-m_2^{p+1}}{k-m_2}-\sum_{0\leq p\leq l}\frac{k^{l-p}}{p+1}\left( 1+\frac{1}{2}+\cdots+\frac{1}{p+1}\right)       k^{p+1}\Bigg{|}\\  
             \end{split}
           \]
           \[
           \begin{split}
           &\leq \sum_{0\leq p\leq l}\frac{k^{l-p}}{p+1}  \Bigg{|} \sum_{1\leq m_2<k}\frac{k^{p+1}-m_2^{p+1}}{k-m_2}- \left(1+\frac{1}{2}+\cdots+\frac{1}{p+1}\right)k^{p+1}     \Bigg{|}    \\
           &\leq k^l+(M+1)\sum_{1\leq p\leq l}\frac{p}{p+1}k^l\\
           &\leq  (M+1)(l+1)k^l.         \\
                      \end{split}
                             \]
                             Since $|A-C|\leq |A-B|+|B-C|$, the lemma is proved.  $\hfill\Box$\\  
   
     \begin{lem}\label{trih}
     For $k>1$, define 
\[N_1(k):=\sum_{\substack{1\leq m_3\leq  m_2<k   \\1\leq m_4\leq m_2<k}}\frac{1}{m_2m_3m_4},\]
we have
\[
N_1(k)=  \frac{1}{3}\left( \sum_{1\leq m<k}\frac{1}{m}\right)^3+\frac{5}{3}\zeta(3)+\mathrm{O}\left( \frac{\mathrm{log}\;k}{k}\right).       
\]
  \end{lem}
  \noindent{\bf Proof:} 
  For $m_3,m_4$, there are three cases: 
  \[ m_3<m_4,\; m_3>m_4,\; m_3=m_4.\]
  Thus 
  \[
\begin{split}
&\;\;\;\;N_1(k)\\
&= \sum_{1\leq m_3<m_4\leq m_2<k}\frac{1}{m_2m_3m_4}  +\sum_{1\leq m_4<m_3\leq m_2<k}\frac{1}{m_2m_3m_4}+\sum_{1\leq m_3\leq m_2<k}\frac{1}{m_2m_3^2}    \\
&= \sum_{1\leq m_3<m_4<m_2<k}\frac{1}{m_2m_3m_4}  +\sum_{1\leq m_4<m_3< m_2<k}\frac{1}{m_2m_3m_4}+\sum_{1\leq m_3<m_2<k}\frac{1}{m_2m_3^2}    \\
&\;\;\;\;+\sum_{1\leq m_3<m_4<k}\frac{1}{m_3m_4^2}+ \sum_{1\leq m_4<m_3<k}\frac{1}{m_4m_3^2}+\sum_{1\leq m_2<k} \frac{1}{m_2^3}    \\
&=2\sum_{1\leq m_4<m_3<m_2<k}\frac{1}{m_2m_3m_4}+\sum_{1\leq m_3<m_2<k}\frac{1}{m_3^2m_2} +2\zeta(1,2)+\zeta(3)+\mathrm{O}\left(\frac{1}{k}\right).\\                              
\end{split}
\]
By the theory of regularization of double shuffle relation in Section $1$, \cite{ikz}, one has
\[
\begin{split}
&\;\;\;\;N_1(k)\\
&= \frac{1}{3}\left( \sum_{1\leq m<k}\frac{1}{m}\right)^3-\zeta(2)  \left( \sum_{1\leq m<k}\frac{1}{m}\right)+\frac{2}{3}\zeta(3) +\zeta(2) \left( \sum_{1\leq m<k}\frac{1}{m}\right) -2\zeta(3)\\
&\;\;\;\;\;\;\;\;\;\;\;\;\;\;\;+2\zeta(1,2)+\zeta(3)+\mathrm{O}\left(\frac{\mathrm{log}\;k}{k}\right)  \\
&=  \frac{1}{3}\left( \sum_{1\leq m<k}\frac{1}{m}\right)^3+\frac{5}{3}\zeta(3)+\mathrm{O}\left( \frac{\mathrm{log}\;k}{k}\right).\\
\end{split}
\]
 $\hfill\Box$\\   
 
 \begin{lem}\label{triha}
\[ \sum_{1\leq q_1<q_2}\frac{1}{q_1q_2}\frac{1}{2^{q_2}}\sum_{1\leq n\leq q_2}\frac{1}{n}=-\frac{\zeta(3)}{8}+\frac{\zeta(2)\,\mathrm{log}\,2}{2}.     \]
 \end{lem}
  \noindent{\bf Proof:} 
  We have
  \[
  \begin{split}
  &\;\;\;\;\sum_{1\leq q_1<q_2}\frac{1}{q_1q_2}\frac{1}{2^{q_2}}\sum_{1\leq n\leq q_2}\frac{1}{n}\\
 & =  \sum_{1\leq q_1<q_2}\frac{1}{q_1q_2}\frac{1}{2^{q_2}}\left(\sum_{1\leq n< q_2}\frac{1}{n} +\frac{1}{q_2}\right)      \\
 &=\left(\sum_{1\leq q_1<n<q_2}+\sum_{1\leq n<q_1<q_2}+\sum_{1\leq q_1=n<q_2}\right)\frac{1}{nq_1q_2}\frac{1}{2^{q_2}} +\sum_{1\leq q_1<q_2}\frac{1}{q_1q_2^2}\frac{1}{2^{q_2}}   \\
 &= 2\sum_{1\leq n<q_1<q_2}\frac{1}{nq_1q_2}\frac{1}{2^{q_2}}+\sum_{1\leq q_1<q_2}\frac{1}{q_1^2q_2}\frac{1}{2^{q_2}}+\sum_{1\leq q_1<q_2}\frac{1}{q_1q_2^2}\frac{1}{2^{q_2}}.   \\
 \end{split}
  \]
  By the Appendix in \cite{zlo}, one has 
  \[
  \begin{split}
  &\sum_{1\leq n<q_1<q_2}\frac{1}{nq_1q_2}\frac{1}{2^{q_2}} =\frac{\mathrm{log}^32}{6},
  \sum_{1\leq q_1<q_2}\frac{1}{q_1^2q_2}\frac{1}{2^{q_2}}=-\frac{\zeta(3)}{4}+\frac{\zeta(2)\,\mathrm{log}\,2}{2}-\frac{\mathrm{log}^32}{6},\\
     &\sum_{1\leq q_1<q_2}\frac{1}{q_1q_2^2}\frac{1}{2^{q_2}}=\frac{\zeta(3)}{8} -\frac{\mathrm{log}^32}{6}.\\
     \end{split} \]
     So we have
     \[
      \sum_{1\leq q_1<q_2}\frac{1}{q_1q_2}\frac{1}{2^{q_2}}\sum_{1\leq n\leq q_2}\frac{1}{n}=-\frac{\zeta(3)}{8}+\frac{\zeta(2)\,\mathrm{log}\,2}{2}.     \]
   
  $\hfill\Box$\\

  \begin{lem}\label{twis}
  For $k\geq 1$, denote by 
  \[
  H(k):=\sum_{1\leq m_i\leq k}\frac{1}{m_2(m_1+m_2)(m_2+m_3)},
  \]
  then
  \[
  H(k)=   \frac{1}{3}\left( \sum_{1\leq m\leq k}\frac{1}{m}\right)^3-\zeta(2)  \left(\sum_{1\leq m\leq k}\frac{1}{m} \right) +\frac{41}{12}\zeta(3)+\mathrm{O}\left( \frac{\mathrm{log}^3\;k}{k}\right). 
 \]
     \end{lem}
     \noindent{\bf Proof:} 
  We have 
 \[
\begin{split}
&\;\;\;\;H(k)   \\
&= \sum_{1\leq m_2\leq k}\frac{1}{m_2} \left(\sum_{1\leq m\leq k}\frac{1}{m+m_2}     \right)^2 \\
&= \sum_{1\leq m_2\leq k}\frac{1}{m_2} \left(\sum_{m_2+1\leq m\leq k+m_2}\frac{1}{m}     \right)^2\\
  \end{split}
           \]
           \[
           \begin{split}
&=  \sum_{1\leq m_2\leq k}\frac{1}{m_2} \left[ \sum_{1\leq m\leq k}\frac{1}{m}+\sum_{1\leq m\leq m_2}\left( \frac{1}{k+m}-\frac{1}{m}    \right)   \right]^2  \\
&=  \sum_{1\leq m_2\leq k}\frac{1}{m_2} \Bigg{[}\left(\sum_{1\leq m\leq k}\frac{1}{m}\right)^2+2\sum_{1\leq m_1\leq k}\frac{1}{m_1}\cdot\sum_{1\leq m_3\leq m_2}\left( \frac{1}{k+m_3}-\frac{1}{m_3}    \right) \\
& \;\;\;\;\;\;\;\;\;\; +\sum_{1\leq m_3,m_4\leq m_2}\left( \frac{1}{k+m_3}-\frac{1}{m_3}    \right) \left( \frac{1}{k+m_4}-\frac{1}{m_4}    \right)  \Bigg{] } \\
&= \left(\sum_{1\leq m\leq k}\frac{1}{m}  \right)^3 +2\sum_{1\leq m\leq k}\frac{1}{m} \cdot \sum_{1\leq m_3\leq  m_2\leq k}\frac{1}{m_2}\left( \frac{1}{k+m_3}-\frac{1}{m_3}   \right)       \\
&\;\;\;\;\;\;\;\;\;\;+\sum_{\substack{  1\leq m_3\leq m_2\leq k,   \\    1\leq m_4\leq m_2\leq k     }}\frac{k^2}{m_2m_3m_4(k+m_3)(k+m_4)}.
\end{split}
\]

Since 
\[
\sum_{\substack{  1\leq m_3\leq m_2=k,   \\    1\leq m_4\leq m_2=k     }}\frac{k^2}{m_2m_3m_4(k+m_3)(k+m_4)}=\mathrm{O}\left( \frac{\mathrm{log}^2k}{k}   \right),
\]

we have
\[
\begin{split}
&\;\;\;\;H(k)\\
&= \left(\sum_{1\leq m\leq k}\frac{1}{m}  \right)^3 +2\sum_{1\leq m\leq k}\frac{1}{m} \cdot \sum_{1\leq m_3\leq  m_2\leq k}\frac{1}{m_2}\left( \frac{1}{k+m_3}-\frac{1}{m_3}   \right)       \\
&\;\;\;\;\;\;\;\;\;\;\;\;\;\;\;\;\;\;\;\;\;\;\;\;\;\;\;\;\;\;\;\;\;+\sum_{\substack{  1\leq m_3\leq m_2< k,   \\    1\leq m_4\leq m_2< k     }}\frac{k^2}{m_2m_3m_4(k+m_3)(k+m_4)}+\mathrm{O}\left( \frac{\mathrm{log}^2k}{k}\right) .  \\
\end{split}
\]
Denote by 
\[
M(k):=\sum_{1\leq m_3\leq  m_2\leq k}\frac{1}{m_2(k+m_3)}  ,
\]
\[
N(k):=\sum_{\substack{  1\leq m_3\leq m_2< k,   \\    1\leq m_4\leq m_2< k     }}\frac{k^2}{m_2m_3m_4(k+m_3)(k+m_4)}.
\]
By the following observation
\[
\sum_{1\leq m_3\leq m_2= k}\frac{1}{m_2(k+m_3)}=\frac{1}{k}\sum_{k+1\leq m\leq 2k}\frac{1}{m}=\mathrm{O}\left( \frac{1}{k}\right),
\]
it follows that
\[
M(k)=\sum_{1\leq m_3\leq  m_2<k}\frac{1}{m_2(k+m_3)}  +\mathrm{O}\left( \frac{1}{k}\right).
\]
So 
\[
\begin{split}
&\;\;\;\;M(k)\\
&=\sum_{1\leq m_2\leq m_3<k}\frac{1}{(k-m_2)(2k-m_3)}+  \mathrm{O}\left( \frac{1}{k}\right)         \\
&=\sum_{1\leq m_2\leq m_3<k}\frac{1}{2k(k-m_2)(1-\frac{m_3}{2k})} +  \mathrm{O}\left( \frac{1}{k}\right)    \\
&=  \sum_{1\leq m_2\leq m_3<k}\frac{1}{2k(k-m_2)}+\sum_{l\geq1}  \sum_{1\leq m_2\leq m_3<k}\frac{1}{k-m_2}\cdot \frac{m_3^l}{(2k)^{l+1}}+  \mathrm{O}\left( \frac{1}{k}\right) \\
&=\frac{1}{2}+\sum_{l\geq 1}\frac{1}{(2k)^{l+1}}\sum_{1\leq m_2\leq m_3<k}\frac{m_3^l}{k-m_2}+\mathrm{O}\left(\frac{1}{k}\right).    \\
\end{split}
\]

By Lemma \ref{pow}, one has
\[
\Bigg{|}\sum_{m_2\leq m_3<k}m_3^l-\left(\frac{k^{l+1}}{l+1}-\frac{k^l}{2}-\frac{m_2^{l+1}}{l+1}+\frac{m_2^l}{2}\right) \Bigg{|}
\leq 2Mlk^{l-1}  .  \tag{4}  
\]
By the above inequation, one has
\[
\Bigg{|}\sum_{1\leq m_2\leq m_3<k}\frac{m_3^l}{k-m_2} -\sum_{1\leq m_2<k} \left( \frac{1}{l+1}\frac{k^{l+1}-m_2^{l+1}}{k-m_2} -\frac{1}{2}\frac{k^l-m_2^l}{k-m_2}   \right)  \Bigg{|}\leq 2Mlk^l.
\]
Thus 
\[
\begin{split}
&\Bigg{|} \sum_{l\geq 1}\frac{1}{(2k)^{l+1}}  \sum_{1\leq m_2\leq m_3<k}\frac{m_3^l}{k-m_2} - \sum_{l\geq 1}\frac{1}{(2k)^{l+1}}   \sum_{1\leq m_2<k} \left( \frac{1}{l+1}\frac{k^{l+1}-m_2^{l+1}}{k-m_2} -\frac{1}{2}\frac{k^l-m_2^l}{k-m_2}   \right)  \Bigg{|}\\
&\leq \sum_{l\geq 1}\frac{1}{(2k)^{l+1}} \cdot 2Mlk^l=\frac{M}{k}\sum_{l\geq 1}\frac{l}{2^l} =\frac{2M}{k} . \\
\end{split}
\]
As a result, $M(k)$ can be reduced to 
\[
M(k)=\frac{1}{2}+\sum_{l\geq 1}\frac{1}{(2k)^{l+1}}   \sum_{1\leq m_2<k} \left( \frac{1}{l+1}\frac{k^{l+1}-m_2^{l+1}}{k-m_2} -\frac{1}{2}\frac{k^l-m_2^l}{k-m_2}   \right)+\mathrm{O}\left(\frac{1}{k}\right).
\]
By Lemma \ref{dpow}, one has
\[
\begin{split}
&\;\;\;\;M(k)\\
&=\frac{1}{2}+\sum_{l\geq 1}\frac{1}{2^{l+1}}  \cdot \frac{1}{l+1}\left(1+\frac{1}{2}+\cdots+\frac{1}{l+1}\right)+\mathrm{O}\left(\frac{1}{k}\right)\\
&=  \sum_{1\leq l^\prime \leq l }\frac{1}{l^\prime l}\frac{1}{2^l}  +\mathrm{O}\left(\frac{1}{k}\right) =\mathrm{Li}_2\left( \frac{1}{2}\right) +\sum_{1\leq l_1<l_2}\frac{1}{l_1l_2}\frac{1}{2^{l_2}}  +\mathrm{O}\left(\frac{1}{k}\right).\\   \end{split}  \tag{5}
\]
 By the formula $ (7.5)$ in \cite{bbb}, we have
 \[
 \mathrm{Li}_2\left(\frac{1}{2}\right)=\sum_{l\geq 1}\frac{1}{l^2}\frac{1}{2^l}=\frac{1}{12}\pi^2-\frac{1}{2}\mathrm{log}^22.
 \]
 For the second term of $M{(k)}$, one has
 \[
 \begin{split}
 &\sum_{1\leq l_1<l_2}\frac{1}{l_1l_2}\frac{1}{2^{l_2}}=\sum_{ l_1,l_2\geq 1}\frac{1}{l_1(l_1+l_2)}\frac{1}{2^{l_1+l_2}} \\
 &=\frac{1}{2}\sum_{l_1,l_2\geq 1}\left( \frac{1}{l_1(l_1+l_2)}+\frac{1}{l_2(l_1+l_2)}   \right)\frac{1}{2^{l_1+l_2}}\\
 &=\frac{1}{2}\sum_{l_1,l_2\geq 1}\frac{1}{l_1l_2}\frac{1}{2^{l_1+l_2} } =\frac{1}{2}\left(\sum_{l\geq 1} \frac{1}{l}\frac{1}{2^l}  \right)^2=\frac{1}{2}\mathrm{log}^22. \\
 \end{split}
 \]
 As a result, \[M(k)=\frac{\pi^2}{12}+\mathrm{O}\left(\frac{1}{k}  \right). \tag{6}\]
 
 We have
 \[
 \begin{split}
 &\;\;\;\;N(k)        \\
 &=\sum_{\substack{1\leq m_3\leq  m_2<k   \\1\leq m_4\leq m_2<k}} \frac{1}{m_2}\left(  \frac{1}{m_3}-\frac{1}{k+m_3}   \right)\left(  \frac{1}{m_4}-\frac{1}{k+m_4}   \right)     \\
 &=\sum_{\substack{1\leq m_3\leq  m_2<k   \\1\leq m_4\leq m_2<k}}\left(   \frac{1}{m_2m_3m_4}-\frac{1}{m_2(k+m_3)m_4}  -\frac{1}{m_2m_3(k+m_4)}+\frac{1}{m_2(k+m_3)(k+m_4)} \right)  \\
 &=N_1(k)-2N_2(k)+N_3(k),
   \end{split}
  \]
  where 
  \[
  N_1(k)=\sum_{\substack{1\leq m_3\leq  m_2<k   \\1\leq m_4\leq m_2<k}}\frac{1}{m_2m_3m_4},
  \]
  \[
  N_2(k)=\sum_{\substack{1\leq m_3\leq  m_2<k   \\1\leq m_4\leq m_2<k}}\frac{1}{m_2(k+m_3)m_4}=\sum_{\substack{1\leq m_3\leq  m_2<k   \\1\leq m_4\leq m_2<k}}  \frac{1}{m_2m_3(k+m_4)}     ,
      \]
      \[
      N_3(k)=\sum_{\substack{1\leq m_3\leq  m_2<k   \\1\leq m_4\leq m_2<k}}\frac{1}{m_2(k+m_3)(k+m_4)}   .
            \]
            
            By Lemma \ref{trih}, one has
            \[
N_1(k)=  \frac{1}{3}\left( \sum_{1\leq m<k}\frac{1}{m}\right)^3+\frac{5}{3}\zeta(3)+\mathrm{O}\left( \frac{\mathrm{log}\;k}{k}\right).       
\]  
One can also use the  method in the calculation of $M(k)$ to  calculate $N_2(k)$ and $N_3(k)$.   Thus
\[
\begin{split}
&\;\;\;\;N_2(k)\\
&=\sum_{\substack{1\leq m_2\leq  m_3<k   \\1\leq m_2\leq m_4<k}} \frac{1}{(k-m_2)(2k-m_3)(k-m_4)}       \\
  \end{split}
           \]
           \[
           \begin{split}
&=   \sum_{\substack{1\leq m_2\leq  m_3<k   \\1\leq m_2\leq m_4<k}}  \frac{1}{(k-m_2)(k-m_4)} \sum_{p\geq 0} \frac{m_3^{p}}{ (2k)^{p+1}}      \\
&= \sum_{p\geq 0} \frac{1}{(2k)^{p+1}} \sum_{\substack{1\leq m_2\leq  m_3<k   \\1\leq m_2\leq m_4<k}}    \frac{m_3^p}{(k-m_2)(k-m_4)}  \\
&=\frac{1}{2k} \sum_{\substack{1\leq m_2\leq  m_3<k   \\1\leq m_2\leq m_4<k}}\frac{1}{(k-m_2)(k-m_4)}    +\sum_{p\geq 1} \frac{1}{(2k)^{p+1}} \sum_{\substack{1\leq m_2\leq  m_3<k   \\1\leq m_2\leq m_4<k}}    \frac{m_3^p}{(k-m_2)(k-m_4)} \\
\end{split}
\]
\[
\begin{split}
&=\frac{1}{2k} \sum_{1\leq m_2\leq m_4<k}\frac{1}{k-m_4}+ \sum_{p\geq 1} \frac{1}{(2k)^{p+1}} \sum_{\substack{1\leq m_2\leq  m_3<k   \\1\leq m_2\leq m_4<k}}    \frac{m_3^p}{(k-m_2)(k-m_4)}      \\
&= \frac{1}{2k}\sum_{1\leq m_4<k}\frac{m_4}{k-m_4}+    \sum_{p\geq 1} \frac{1}{(2k)^{p+1}} \sum_{\substack{1\leq m_2\leq  m_3<k   \\1\leq m_2\leq m_4<k}}    \frac{m_3^p}{(k-m_2)(k-m_4)}      \\
  &=  \frac{1}{2k}\sum_{1\leq m_4<k}\left(-1+\frac{k}{k-m_4}\right)+    \sum_{p\geq 1} \frac{1}{(2k)^{p+1}} \sum_{\substack{1\leq m_2\leq  m_3<k   \\1\leq m_2\leq m_4<k}}    \frac{m_3^p}{(k-m_2)(k-m_4)}      \\
  &=-\frac{1}{2}+\frac{1}{2}\sum_{1\leq m\leq k}\frac{1}{m}+\mathrm{O}\left(\frac{1}{k}\right)+  \sum_{p\geq 1} \frac{1}{(2k)^{p+1}} \sum_{\substack{1\leq m_2\leq  m_3<k   \\1\leq m_2\leq m_4<k}}    \frac{m_3^p}{(k-m_2)(k-m_4)}      \\
  &=-\frac{1}{2}+\frac{1}{2}\sum_{1\leq m\leq k}\frac{1}{m}+  \sum_{p\geq 1} \frac{1}{(2k)^{p+1}} \sum_{\substack{1\leq m_2\leq  m_3<k   \\1\leq m_2\leq m_4<k}}    \frac{m_3^p}{(k-m_2)(k-m_4)}+\mathrm{O}\left(\frac{1}{k}\right)  .    \\
  \end{split}
\]
By the formula $(4)$, it follows that
\[
\begin{split}
&\Bigg{|} \sum_{\substack{1\leq m_2\leq  m_3<k   \\1\leq m_2\leq m_4<k}}    \frac{m_3^p}{(k-m_2)(k-m_4)}-\sum_{1\leq m_2\leq m_4<k}\frac{1}{(k-m_2)(k-m_4)}\left(\frac{k^{p+1}-m_2^{p+1}}{p+1}-\frac{k^p-m_2^p}{2}     \right)\\
&\leq \sum_{1\leq m_2\leq m_4<k}\frac{2Mpk^{p-1}}{(k-m_2)(k-m_4)}=\sum_{1\leq m_4\leq m_2<k}\frac{2Mpk^{p-1}}{m_2m_4}.\\
\end{split}
\]

So we have
\[
\begin{split}
&\Bigg{|}\sum_{p\geq 1}\frac{1}{(2k)^{p+1}}\sum_{\substack{1\leq m_2\leq  m_3<k   \\1\leq m_2\leq m_4<k}}   \frac{m_3^p}{(k-m_2)(k-m_4)}\\
&\;\;\;\;- \sum_{p\geq 1}\frac{1}{(2k)^{p+1}}     \sum_{1\leq m_2\leq m_4<k}\frac{1}{(k-m_2)(k-m_4)}\left(\frac{k^{p+1}-m_2^{p+1}}{p+1}-\frac{k^p-m_2^p}{2}       \right)\Bigg{|}\\
&\leq \sum_{p\geq 1}\frac{1}{(2k)^{p+1}} \sum_{1\leq m_4\leq m_2<k}\frac{2Mpk^{p-1}}{m_2m_4}=\mathrm{O}\left(\frac{\mathrm{log}^2k}{k^2}\right).\\
\end{split}
\]
From the above inequation, one has 
\[
\begin{split}
&\;\;\;\;N_2(k)\\
&=  -\frac{1}{2}+\frac{1}{2}\sum_{1\leq m\leq k}\frac{1}{m}+   \\
&+ \sum_{p\geq 1}\frac{1}{(2k)^{p+1}}     \sum_{1\leq m_2\leq m_4<k}\frac{1}{(k-m_2)(k-m_4)}\left(\frac{k^{p+1}-m_2^{p+1}}{p+1}-\frac{k^p-m_2^p}{2}       \right)   + \mathrm{O}\left(\frac{1}{k}\right).   \\
\end{split}
\]

By Lemma \ref{ddpow}, we have
\[
\begin{split}
&\;\;\;\;N_2(k)\\
&= -\frac{1}{2}+\frac{1}{2}\sum_{1\leq m\leq k}\frac{1}{m}\\
&+\sum_{p\geq1}\frac{1}{2^{p+1}(p+1)}\left(\sum_{0\leq p_1\leq p}\frac{1}{p_1+1}\sum_{1\leq m<k}\frac{1}{m}-  \sum_{0\leq p_1\leq p}\frac{1}{p_1+1}(1+\frac{1}{2}+\cdots+\frac{1}{p_1+1}                  )            \right)\\
& +\mathrm{O}\left( \frac{\mathrm{log}\;k}{k}\right)\\
&=\sum_{p\geq0}\frac{1}{2^{p+1}(p+1)}\left(\sum_{0\leq p_1\leq p}\frac{1}{p_1+1}\sum_{1\leq m<k}\frac{1}{m}-  \sum_{0\leq p_1\leq p}\frac{  1+\frac{1}{2}+\cdots+\frac{1}{p_1+1}                     }{p_1+1}          \right) +\mathrm{O}\left( \frac{\mathrm{log}\;k}{k}\right)\\
&= \sum_{1\leq p_1\leq p}\frac{1}{p_1p}\frac{1}{2^p} \left(\sum_{1\leq m<k}\frac{1}{m} \right) -\sum_{1\leq p_2\leq p_1\leq p}\frac{1}{p_2p_1p}\frac{1}{2^p} +\mathrm{O}\left(\frac{\mathrm{\log}\;k }{k}\right).             \\
 \end{split}
\]
By the formulas $(5)$ and $(6)$, we have
\[
N_2(k)=\frac{1}{2}\zeta(2) \left(\sum_{1\leq m<k}\frac{1}{m} \right)-\sum_{1\leq p_2\leq p_1\leq p}\frac{1}{p_2p_1p}\frac{1}{2^p}+\mathrm{O}\left( \frac{\mathrm{log}\;k}{k}\right).
\]
By Theorem $3$ in \cite{zlo}, one has 
\[
N_2(k)=\frac{1}{2}\zeta(2) \left(\sum_{1\leq m<k}\frac{1}{m} \right)-\frac{3}{4}\zeta(3)+\mathrm{O}\left( \frac{\mathrm{log}\;k}{k}\right).
\]

For $N_3(k)$, we have
\[
\begin{split}
  &\;\;\;\;N_3(k)\\
  &=\sum_{\substack{1\leq m_3\leq  m_2<k   \\1\leq m_4\leq m_2<k}}\frac{1}{m_2(k+m_3)(k+m_4)}   \\
  &= \sum_{\substack{1\leq m_2\leq  m_3<k   \\1\leq m_2\leq m_4<k}} \frac{1}{(k-m_2)(2k-m_3)(2k-m_4)}     \\
  &=\sum_{\substack{1\leq m_2\leq  m_3<k   \\1\leq m_2\leq m_4<k}}\frac{1}{k-m_2}\sum_{q_1,q_2\geq 0}\frac{m_3^{q_1}m_4^{q_2}}{(2k)^{q_1+q_2+2}}.\\
    \end{split}
  \]
  By the formula $(4)$ (it also holds for $l=0$), one has 
  \[
  \begin{split}
  &\;\;\;\;\Bigg{|}N_3(k)-\sum_{\substack{1\leq m_2<k\\ q_1,q_2\geq 0}} \frac{  \left(\frac{k^{q_1+1}-m_2^{q_1+1}}{q_1+1}-\frac{k^{q_1}-m_2^{q_1}}{2}      \right)\left(\frac{k^{q_2+1}-m_2^{q_2+1}}{q_2+1}-\frac{k^{q_2}-m_2^{q_2}}{2}      \right)     }{(k-m_2)(2k)^{q_1+q_2+2}} \Bigg{|}  \\
  &\leq  \Bigg{|} N_3(k)- \sum_{\substack{1\leq m_2\leq m_4<k\\ q_1,q_2\geq 0}} \frac{\left(\frac{k^{q_1+1}-m_2^{q_1+1}}{q_1+1}-\frac{k^{q_1}-m_2^{q_1}}{2}      \right) m_4^{q_2}}{(k-m_2)(2k)^{q_1+q_2+2}} \Bigg{|}\\
  &+\Bigg{|}  \sum_{\substack{1\leq m_2<k\\ q_1,q_2\geq 0}} \frac{\left(\frac{k^{q_1+1}-m_2^{q_1+1}}{q_1+1}-\frac{k^{q_1}-m_2^{q_1}}{2}      \right) \left(\sum_{m_2\leq m_4<k}m_4^{q_2}-( \frac{k^{q_2+1}-m_2^{q_2+1}}{q_2+1}-\frac{k^{q_2}-m_2^{q_2}}{2}             ) \right)     }{(k-m_2)(2k)^{q_1+q_2+2}}      \Bigg{|}\\
  &\leq  \sum_{\substack{1\leq m_2\leq m_4<k\\ q_1,q_2\geq 0}}\frac{2Mq_1k^{q_1-1}m_4^{q_2}}{(k-m_2)(2k)^{q_1+q_2+2}}+ \sum_{\substack{1\leq m_2<k\\ q_1,q_2\geq 0}} \frac{\left(\frac{k^{q_1+1}-m_2^{q_1+1}}{q_1+1}-\frac{k^{q_1}-m_2^{q_1}}{2}      \right) \cdot 2Mq_2k^{q_2-1}    }{(k-m_2)(2k)^{q_1+q_2+2}}     \\
  &\leq \sum_{\substack{1\leq m_2<k\\q_1,q_2\geq0}}\frac{2Mq_1k^{q_1+q_2}}{(k-m_2)(2k)^{q_1+q_2+2}}+\sum_{\substack{1\leq m_2<k\\ q_1,q_2\geq 0}} \frac{2Mq_2k^{q_1+q_2}}{(k-m_2)(2k)^{q_1+q_2+2}}  \\
 & =\mathrm{O}\left(\frac{\mathrm{log}\;k}{k^2}\right).\\ \end{split}
  \]
  
  By Lemma \ref{bipow}, it follows that
  \[
  \begin{split}
  &\;\;\;\;N_3(k)\\
  & = \sum_{q_1,q_2\geq 0}\frac{  \sum_{1\leq n_1\leq q_1+1}\frac{1}{n_1}+ \sum_{1\leq n_2\leq q_2+1}\frac{1}{n_2}-\sum_{1\leq n_3\leq q_1+q_2+2}\frac{1}{n_3}       }{2^{q_1+q_2+2}(q_1+1)(q_2+1)} +\mathrm{O}\left( \frac{1}{k}\right)\\
  &= \sum_{q_1,q_2\geq 1}\frac{1}{2^{q_1+q_2}q_1q_2} \left(\sum_{1\leq n_1\leq q_1}\frac{1}{n_1}+ \sum_{1\leq n_2\leq q_2}\frac{1}{n_2}-\sum_{1\leq n_3\leq q_1+q_2}\frac{1}{n_3}\right) +\mathrm{O}\left( \frac{1}{k}\right)   \\
  &=2\sum_{\substack{q_1\geq n_1\geq 1\\q_2\geq 1}}\frac{1}{n_1q_1q_2}\frac{1}{2^{q_1+q_2}} -\sum_{q_1,q_2\geq1}\frac{1}{2^{q_1+q_2}q_1q_2}\sum_{1\leq n_3\leq q_1+q_2}\frac{1}{n_3} + \mathrm{O}\left( \frac{1}{k}\right)   \\
  &=2\sum_{q_1\geq n_1\geq 1}\frac{1}{n_1q_12^{q_1}}\sum_{q_2\geq 1}\frac{1}{q_2 2^{q_2}}-\sum_{q_1,q_2\geq1}\frac{ \frac{1}{q_1(q_1+q_2)}+ \frac{1}{q_2(q_1+q_2)} }{2^{q_1+q_2}}\sum_{1\leq n_3\leq q_1+q_2}\frac{1}{n_3} + \mathrm{O}\left( \frac{1}{k}\right)\\  
   &=\zeta(2)\,\mathrm{log}\;2-2\sum_{q_1,q_2\geq1}\frac{ {1} }{ {q_1(q_1+q_2)}2^{q_1+q_2}}\sum_{1\leq n_3\leq q_1+q_2}\frac{1}{n_3} + \mathrm{O}\left( \frac{1}{k}\right)\\
      \end{split}
  \]
  \[
  \begin{split}   &=\zeta(2)\,\mathrm{log}\;2-2\sum_{1\leq q_1<q_2} \frac{1}{q_1q_22^{q_2}}\sum_{1\leq n_3\leq q_2}\frac{1}{n_3}  + \mathrm{O}\left( \frac{1}{k}\right)\\
  &=\zeta(2)\,\mathrm{log}\;2-2( -\frac{\zeta(3)}{8}+\frac{\zeta(2)\,\mathrm{log}\,2}{2})+ \mathrm{O}\left( \frac{1}{k}\right)\\
  &=\frac{\zeta(3)}{4}+\mathrm{O}\left( \frac{1}{k}\right)
          \end{split}
  \]
   
  In conclusion, we have 
  \[
  \begin{split}
  &\;\;\;\;N(k)\\
  &=N_1(k)-2N_2(k)+N_3(k)\\
  &= \frac{1}{3}\left( \sum_{1\leq m<k}\frac{1}{m}\right)^3+\frac{5}{3}\zeta(3)-2\left(     \frac{1}{2}\zeta(2) \left(\sum_{1\leq m<k}\frac{1}{m} \right)-\frac{3}{4}\zeta(3)     \right)  +\frac{\zeta(3)}{4}+\mathrm{O}\left( \frac{\mathrm{log}\;k}{k}\right)  \\
  &=  \frac{1}{3}\left( \sum_{1\leq m<k}\frac{1}{m}\right)^3-\zeta(2)  \left(\sum_{1\leq m<k}\frac{1}{m} \right) +\frac{41}{12}\zeta(3)+\mathrm{O}\left( \frac{\mathrm{log}\;k}{k}\right),    \\
    \end{split}
  \]
  \[
 \begin{split}
 &\;\;\;\;H(k)\\
 &=\left( \sum_{1\leq m\leq k}\frac{1}{m}\right)^3+2\sum_{1\leq m\leq k}\frac{1}{m}\cdot \left(M(k)-\sum_{1\leq m_3\leq m_2\leq k}   \frac{1}{m_2m_3}     \right)  +N(k)          \\
 &= \left( \sum_{1\leq m\leq k}\frac{1}{m}\right)^3+2\sum_{1\leq m\leq k}\frac{1}{m}\cdot \left(M(k)-\frac{1}{2}\left( \sum_{1\leq m\leq k}\frac{1}{m}  \right)^2-\frac{1}{2}\sum_{1\leq m\leq k}\frac{1}{m^2}  \right)  +N(k) \\
 &=[2M(k)-\zeta(2)]\left( \sum_{1\leq m\leq k}\frac{1}{m}\right)+N(k)+\mathrm{O}\left(\frac{\mathrm{log}\;k}{k}\right)   \\
 &=N(k)+\mathrm{O}\left(\frac{\mathrm{log}\;k}{k}\right)   \\
 &=   \frac{1}{3}\left( \sum_{1\leq m\leq k}\frac{1}{m}\right)^3-\zeta(2)  \left(\sum_{1\leq m\leq k}\frac{1}{m} \right) +\frac{41}{12}\zeta(3)+\mathrm{O}\left( \frac{\mathrm{log}^3\;k}{k}\right).
    \end{split}
 \]
   $\hfill\Box$\\

  \begin{Thm}\label{sum}
  For $k\geq 1$, we have 
  \[
 \sum_{1\leq m_i\leq k}\frac{1}{(m_1+m_2)(m_2+m_3)(m_3+m_1)   }=\frac{3}{2}\zeta(2) \left(\sum_{1\leq m\leq k}\frac{1}{m} \right)-\frac{41}{8}\zeta(3)+ \mathrm{O}\left( \frac{\mathrm{log}^3\;k}{k}\right) .   \]
  \end{Thm}
    \noindent{\bf Proof:} We have
    \[
    \begin{split}
    &\;\;\;\;\sum_{1\leq m_i\leq k}\frac{1}{(m_1+m_2)(m_2+m_3)(m_3+m_1)}\\
    &=\frac{1}{2} \sum_{1\leq m_i\leq k}  \frac{(m_1+m_2)+(m_2+m_3)+(m_3+m_1)}{(m_1+m_2)(m_2+m_3)(m_3+m_1)(m_1+m_2+m_3)}   \\
        \end{split}
    \]
    \[
    \begin{split}    &=\frac{1}{2} \sum_{1\leq m_i\leq k} \Bigg{(} \frac{1}{(m_2+m_3)(m_3+m_1)(m_1+m_2+m_3)}\\
    &+\frac{1}{(m_1+m_2)(m_3+m_1)(m_1+m_2+m_3)}+ \frac{1}{(m_1+m_2)(m_2+m_3)(m_1+m_2+m_3)}\Bigg{)}  \\
    &=\frac{3}{2}\sum_{1\leq m_i\leq k}  \frac{1}{(m_1+m_2)(m_2+m_3)(m_1+m_2+m_3)}\\
    &=\frac{3}{2}\sum_{1\leq m_i\leq k}\frac{1}{m_2}\frac{ (m_1+m_2)+(m_2+m_3)-(m_1+m_2+m_3)}{(m_1+m_2)(m_2+m_3)(m_1+m_2+m_3) }\\
    &=\frac{3}{2}\sum_{1\leq m_i\leq k} \Bigg{(} \frac{1}{m_2(m_2+m_3)(m_1+m_2+m_3)}+\frac{1}{m_2(m_1+m_2)(m_1+m_2+m_3)}\\
    &\;\;\;\;\;\;\;\;\;\;\;\;\;\;\;-\frac{1}{m_2(m_1+m_2)(m_2+m_3)}\Bigg{)}\\
    &=3\sum_{1\leq m_i\leq k}\frac{1}{m_1(m_1+m_2)(m_1+m_2+m_3)}-\frac{3}{2}\sum_{1\leq m_i\leq k}\frac{1}{m_2(m_1+m_2)(m_2+m_3)}.\\
             \end{split}
        \]
        
       From the obvious identity
       \[
       \frac{1}{m_1m_2}=\frac{1}{m_1(m_1+m_2)}+\frac{1}{m_2(m_1+m_2)},
       \]
       one can check that 
       \[
       \frac{1}{m_1m_2m_3}=\sum_{\sigma\in S_3} \frac{1}{m_{\sigma(1)}(m_{\sigma(1)}+m_{\sigma(2)})(m_{\sigma(1)}+m_{\sigma(2)}+m_{\sigma(3)} )}, \tag{7}     \]
       where $S_3 $ is the group of permutations on the set $\{1,2,3\}$.
       So 
       \[
       \begin{split}
       &\;\;\;\;\left( \sum_{1\leq m\leq k}\frac{1}{m} \right)^3=  \sum_{1\leq m_i\leq k}\frac{1}{m_1m_2m_3}\\
        &=\sum_{1\leq m_i\leq k}\sum_{\sigma\in S_3} \frac{1}{m_{\sigma(1)}(m_{\sigma(1)}+m_{\sigma(2)})(m_{\sigma(1)}+m_{\sigma(2)}+m_{\sigma(3)} )}   \\
       &=6 \sum_{1\leq m_i\leq k} \frac{1}{m_{1}(m_{1}+m_{2})(m_{1}+m_{2}+m_{3} )} .   \\
       \end{split}
       \]
       From the above formula and Lemma \ref{twis}, it follows that
       \[
       \begin{split}
       &\;\;\;\; \sum_{1\leq m_i\leq k}\frac{1}{(m_1+m_2)(m_2+m_3)(m_3+m_1)   }\\
       & =\frac{1}{2} \left( \sum_{1\leq m\leq k}\frac{1}{m} \right)^3 -\frac{3}{2}\sum_{1\leq m_i\leq k}\frac{1}{m_2(m_1+m_2)(m_2+m_3)}\\
           \end{split}
    \]
    \[
    \begin{split}       &=\frac{1}{2} \left( \sum_{1\leq m\leq k}\frac{1}{m} \right)^3-\frac{3}{2}H(k)\\
       &=\frac{1}{2}\left(\sum_{1\leq m\leq k}\frac{1}{m}\right)^3-
       \frac{3}{2}\left( \frac{1}{3}\left( \sum_{1\leq m\leq k}\frac{1}{m}\right)^3-\zeta(2)  \left(\sum_{1\leq m\leq k}\frac{1}{m} \right) +\frac{41}{12}\zeta(3)\right)
       + \mathrm{O}\left( \frac{\mathrm{log}^3\;k}{k}\right)   \\  
       &=\frac{3}{2}\zeta(2) \left(\sum_{1\leq m\leq k}\frac{1}{m} \right)-\frac{41}{8}\zeta(3)+ \mathrm{O}\left( \frac{\mathrm{log}^3\;k}{k}\right) . \\   \end{split}
       \]
           $\hfill\Box$\\   
           
           \subsection{The symmetric triple integral}
     
    For the triple integral, we will need the following lemmas.
      \begin{lem}\label{dhtri}
      (i) For $n_1,n_2\geq 1, n_1\neq n_2$, 
  \[
  \begin{split}
  &\int^{k+1}_1\left(\frac{1}{t^{n_1}}-\frac{t^{n_1}}{(k+1)^{n_1}}\right)\left(\frac{1}{t^{n_2}}-\frac{t^{n_2}}{(k+1)^{n_2}}\right) \frac{dt}{t}=\\
  &\;\;\;\;\;\;\;\frac{2}{n_1+n_2}\left(1-\frac{1}{(k+1)^{n_1+n_2}} \right) +\frac{2}{n_1-n_2}\left( \frac{1}{(k+1)^{n_1}}-\frac{1}{(k+1)^{n_2}}      \right) , \\
  \end{split} \]
  (ii) For $n\geq 1$, 
  \[
  \int^{k+1}_1\left(\frac{1}{t^{n}}-\frac{t^{n}}{(k+1)^{n}}\right)\left(\frac{1}{t^{n}}-\frac{t^{n}}{(k+1)^{n}}\right) \frac{dt}{t}= \frac{1}{n}\left( 1-\frac{1}{(k+1)^{2n}}     \right)-\frac{2\mathrm{log}\,(k+1)}{(k+1)^n}.  \]   \end{lem}
  \noindent{\bf Proof:} 
  One can check $(i)$ and $(ii)$ directly. $\hfill\Box$\\

     \begin{lem}\label{htri}
  \[
    \int_1^{k+1}\mathrm{log}^2\left(\frac{k+1+y}{1+y} \right) \frac{dy}{y} =\frac{1}{3}\mathrm{log}^3(k+1)-\zeta(2)\mathrm{log}(k+1)+\frac{7}{2}\zeta(3) +\mathrm{O}\left(\frac{\mathrm{log}\,k}{k}\right) .  \]
   \end{lem}
  \noindent{\bf Proof:}  By letting $y=\frac{1}{t}$, one has 
  \[
  \begin{split}
  &\;\;\;\;\int_1^{k+1}\mathrm{log}^2\left(\frac{k+1+y}{1+y} \right) \frac{dy}{y} \\
  &=\int_{\frac{1}{k+1}}^1\mathrm{log}^2\left(\frac{(k+1)t+1}{1+t} \right) \frac{dt}{t}\\
  &=\int_{\frac{1}{k+1}}^1 \left[ \mathrm{log}[(k+1)t]+\mathrm{log}\left(1+\frac{1}{(k+1)t}\right) -\mathrm{log}(1+t)                  \right] ^2                                  \frac{dt}{t}\\
  &=\int_{\frac{1}{k+1}}^1 \left[ \mathrm{log}[(k+1)t]  +\sum_{n\geq1}\frac{(-1)^{n-1}}{n}\left(\frac{1}{(k+1)^nt^n}-t^n   \right)             \right] ^2  \frac{dt}{t}\\
  &=\int_{\frac{1}{k+1}}^1 \ \mathrm{log}^2[(k+1)t]\frac{dt}{t}+2 \int_{\frac{1}{k+1}}^1 \mathrm{log}[(k+1)t]  \sum_{n\geq1}\frac{(-1)^{n-1}}{n}\left(\frac{1}{(k+1)^nt^n}-t^n   \right) \frac{dt}{t}     \\
    \end{split}
  \]
  \[
  \begin{split}   &+\int_{\frac{1}{k+1}}^1 \sum_{n_1,n_2\geq 1}\frac{(-1)^{n_1+n_2}}{n_1n_2}\left(\frac{1}{(k+1)^{n_1}t^{n_1}}-t^{n _1}  \right) \left(\frac{1}{(k+1)^{n_2}t^{n_2}}-t^{n _2}  \right)     \frac{dt}{t}    \\
  &=\int_{1}^{k+1} \ \mathrm{log}^2t\,\frac{dt}{t}+2 \int_{1}^{k+1} \mathrm{log}\,t \sum_{n\geq1}\frac{(-1)^{n-1}}{n}\left(\frac{1}{t^n}-\frac{t^n}{(k+1)^n}   \right) \frac{dt}{t}     \\
  &\;\;\;\;\;\;\;\;\;+\int_{1}^{k+1} \sum_{n_1,n_2\geq 1}\frac{(-1)^{n_1+n_2}}{n_1n_2}\left(\frac{1}{t^{n_1}}-\frac{t^{n _1}}{(k+1)^{n_1}}  \right) \left(\frac{1}{t^{n_2}}-\frac{t^{n _2}}{(k+1)^{n_2}}  \right)     \frac{dt}{t} .   \\  \end{split}
  \]
  By the following  formulas
  \[
  \int t^{n-1}\mathrm{log}\,tdt=\frac{t^n}{n}\mathrm{log}\,t-\frac{t^n}{n^2}+C,\;n\neq 0,
  \]
  \[
  \int \frac{\mathrm{log}^2t}{t}dt=\frac{1}{3}\mathrm{log}^3t+C,
  \]
   we have
   \[
   \begin{split}
   &\;\;\;\;\int_1^{k+1}\left(\frac{1}{t^n}-\frac{t^n}{(k+1)^n}\right)\frac{\mathrm{log}\,t\,dt}{t}\\
   &=  \left(-\frac{t^{-n}}{n}\mathrm{log}\,t-\frac{t^{-n}}{n^2}\right)\Bigg{|}^{t=k+1}_{t=1} -\frac{1}{(k+1)^n}\left(\frac{t^n}{n}\mathrm{log}\,t-\frac{t^n}{n^2}      \right)\Bigg{|}^{t=k+1}_{t=1}   \\
   &=-\frac{\mathrm{\log}(k+1)}{n}\left(1+\frac{1}{(k+1)^n}\right)+\frac{2}{n^2}\left(1-\frac{1}{(k+1)^n}\right).\\
   \end{split}
   \]
   Thus we have
     \[
     \begin{split}
    &\;\;\;\;\int_1^{k+1}\mathrm{log}^2\left(\frac{k+1+y}{1+y} \right) \frac{dy}{y} \\
   & =  \frac{1}{3}\mathrm{log}^3(k+1)+2\sum_{n\geq 1}\frac{(-1)^{n-1}}{n} \left[  -\frac{\mathrm{\log}(k+1)}{n}\left(1+\frac{1}{(k+1)^n}\right)+\frac{2}{n^2}\left(1-\frac{1}{(k+1)^n}\right)           \right]   \\
   &\;\;\;\;\;\;\;\;\;+ \sum_{n_1,n_2\geq 1} \frac{(-1)^{n_1+n_2}}{n_1n_2}  \int_{1}^{k+1} \left(\frac{1}{t^{n_1}}-\frac{t^{n _1}}{(k+1)^{n_1}}  \right) \left(\frac{1}{t^{n_2}}-\frac{t^{n _2}}{(k+1)^{n_2}}  \right)     \frac{dt}{t} \\
     & =  \frac{1}{3}\mathrm{log}^3(k+1)+2\sum_{n\geq 1}\left[\frac{(-1)^{n}}{n^2}\mathrm{\log}(k+1)+\frac{2(-1)^{n-1}}{n^3}\right]+\mathrm{O}\left(\frac{\mathrm{log}\,k}{k}\right)   \\
   &\;\;\;\;\;\;\;\;\;+ \sum_{n_1,n_2\geq 1} \frac{(-1)^{n_1+n_2}}{n_1n_2}  \int_{1}^{k+1} \left(\frac{1}{t^{n_1}}-\frac{t^{n _1}}{(k+1)^{n_1}}  \right) \left(\frac{1}{t^{n_2}}-\frac{t^{n _2}}{(k+1)^{n_2}}  \right)     \frac{dt}{t} \\
     &=\frac{1}{3}\mathrm{log}^3(k+1)-\zeta(2)\mathrm{log}(k+1)+3\zeta(3) +\mathrm{O}\left(\frac{\mathrm{log}\,k}{k}\right)   \\
 &\;\;\;\;\;\;\;\;\;+ \sum_{n_1,n_2\geq 1} \frac{(-1)^{n_1+n_2}}{n_1n_2}  \int_{1}^{k+1} \left(\frac{1}{t^{n_1}}-\frac{t^{n _1}}{(k+1)^{n_1}}  \right) \left(\frac{1}{t^{n_2}}-\frac{t^{n _2}}{(k+1)^{n_2}}  \right)     \frac{dt}{t} \\
 &=\frac{1}{3}\mathrm{log}^3(k+1)-\zeta(2)\mathrm{log}(k+1)+{3}\zeta(3) +\mathrm{O}\left(\frac{\mathrm{log}\,k}{k}\right)   \\
   \end{split}
  \]
  \[
  \begin{split} 
 &\;\;\;\;\;\;\;\;\;+ 2\sum_{1\leq n_1<n_2} \frac{(-1)^{n_1+n_2}}{n_1n_2}  \int_{1}^{k+1} \left(\frac{1}{t^{n_1}}-\frac{t^{n _1}}{(k+1)^{n_1}}  \right) \left(\frac{1}{t^{n_2}}-\frac{t^{n _2}}{(k+1)^{n_2}}  \right)     \frac{dt}{t} \\ 
  &\;\;\;\;\;\;\;\;\;\;\;\;\;\;\;\;\;\;\;\;\;\;+ \sum_{n\geq 1} \frac{1}{n^2}  \int_{1}^{k+1} \left(\frac{1}{t^{n}}-\frac{t^{n }}{(k+1)^{n}}  \right) \left(\frac{1}{t^{n}}-\frac{t^{n }}{(k+1)^{n}}  \right)     \frac{dt}{t}. \\      
       \end{split}
   \]
  By Lemma \ref{dhtri}, we have
  \[
    \begin{split}
    &\;\;\;\;\int_1^{k+1}\mathrm{log}^2\left(\frac{k+1+y}{1+y} \right) \frac{dy}{y} \\
   &=\frac{1}{3}\mathrm{log}^3(k+1)-\zeta(2)\mathrm{log}(k+1)+{3}\zeta(3) +\mathrm{O}\left(\frac{\mathrm{log}\,k}{k}\right)   \\
   &+ \sum_{1\leq n_1<n_2} \frac{4(-1)^{n_1+n_2}}{n_1n_2} \left[   \frac{1}{n_1+n_2}\left(1-\frac{1}{(k+1)^{n_1+n_2}} \right) +\frac{1}{n_1-n_2}\left( \frac{1}{(k+1)^{n_1}}-\frac{1}{(k+1)^{n_2}}      \right) \right] , \\
   &+\sum_{n\geq 1}\frac{1}{n^2}\left[  \frac{1}{n}\left( 1-\frac{1}{(k+1)^{2n}}     \right)-\frac{2\mathrm{log}\,(k+1)}{(k+1)^n}      \right]\\
   &=\frac{1}{3}\mathrm{log}^3(k+1)-\zeta(2)\mathrm{log}(k+1)+{3}\zeta(3) +\mathrm{O}\left(\frac{\mathrm{log}\,k}{k}\right)   
   +4\sum_{1\leq n_1<n_2}\frac{(-1)^{n_1+n_2}}{n_1n_2(n_1+n_2)}\\
   &+ \mathrm{O}\left(\frac{1}{k}\right)  +\zeta(3)+\mathrm{O}\left(\frac{\mathrm{log}\,k}{k}\right).\\
    \end{split}
     \]
     By the Appendix in \cite{zlo}, we have
     \[
     \begin{split}
     &\;\;\;\;\sum_{1\leq n_1<n_2}\frac{(-1)^{n_1+n_2}}{n_1n_2(n_1+n_2)}\\
     &=\frac{1}{2}\sum_{n_1,n_2\geq 1}\frac{(-1)^{n_1+n_2}}{n_1n_2(n_1+n_2) } -\frac{1}{4}\zeta(3)  \\
     &=\frac{1}{2}\sum_{n_1,n_2\geq 1}\left( \frac{1}{n_1(n_1+n_2)}+\frac{1}{n_2(n_1+n_2)}     \right)\frac{(-1)^{n_1+n_2}}{(n_1+n_2) }  -\frac{1}{4}\zeta(3)\\
     &= \sum_{1\leq n_1<n_2}\frac{(-1)^{n_2}}{n_1n_2^2}-\frac{1}{4}\zeta(3) =\frac{1}{8}\zeta(3)-\frac{1}{4}\zeta(3)=-\frac{1}{8}\zeta(3).
     \end{split}
           \]
           As a result, we have 
            \[
    \begin{split}
    &\int_1^{k+1}\mathrm{log}^2\left(\frac{k+1+y}{1+y} \right) \frac{dy}{y} =\frac{1}{3}\mathrm{log}^3(k+1)-\zeta(2)\mathrm{log}(k+1)+\frac{7}{2}\zeta(3) +\mathrm{O}\left(\frac{\mathrm{log}\,k}{k}\right).  \\
    \end{split}
    \]

 $\hfill\Box$\\

 \begin{Thm}\label{tri}
  \[\mathop{\iiint}_{[1,k+1]^3}\frac{dxdydz}{(x+y)(y+z)(z+x)}= \frac{3}{2}\zeta(2)\mathrm{log}\,(k+1)-\frac{21}{4}\zeta(3)+\mathrm{O}\left(\frac{\mathrm{log}\,k}{k}\right) .  \]
   \end{Thm}
  \noindent{\bf Proof:} 
 Since 
 \[
 \begin{split}
 &\;\;\;\;\frac{1}{(x+y)(y+z)(z+x)}\\
 &= \frac{1}{2}\frac{(x+y)+(y+z)+(z+x)}{(x+y)(y+z)(z+x)(x+y+z)}       \\
 &=\frac{1}{2}\frac{1}{(y+z)(z+x)(x+y+z)}+\frac{1}{2}\frac{1}{(x+y)(z+x)(x+y+z)}+\frac{1}{2}\frac{1}{(x+y)(y+z)(x+y+z)},\\
 \end{split}
 \]
  we have
  \[\mathop{\iiint}_{[1,k+1]^3}\frac{dxdydz}{(x+y)(y+z)(z+x)}=\frac{3}{2} \mathop{\iiint}_{[1,k+1]^3}  \frac{dxdydz}{(x+y)(y+z)(x+y+z)}.  \]
  
  From the identity 
  \[
  \begin{split}
 &\;\;\;\; \frac{1}{(x+y)(y+z)(x+y+z)}\\
 &=\frac{(x+y)+(y+z)-(x+y+z)}{y(x+y)(y+z)(x+y+z)}\\
 &=\frac{1}{y(y+z)(x+y+z)}+\frac{1}{y(x+y)(x+y+z)}-\frac{1}{y(x+y)(y+z)}\\
 \end{split}
  \]
  and the formula $(7)$,  it follows that 
    \[
   \begin{split}
   &\;\;\;\; \mathop{\iiint}_{[1,k+1]^3}\frac{dxdydz}{(x+y)(y+z)(z+x)}\\
   &= \frac{1}{2} \mathop{\iiint}_{[1,k+1]^3}\frac{dxdydz}{xyz}   -\frac{3}{2}    \mathop{\iiint}_{[1,k+1]^3}\frac{dxdydz}{y(x+y)(y+z)}         \\
   &=\frac{1}{2}\mathrm{log}^3(k+1)-\frac{3}{2} \mathop{\iiint}_{[1,k+1]^3}\frac{dxdydz}{y(x+y)(y+z)}     \\
   &=\frac{1}{2}\mathrm{log}^3(k+1)
   -\frac{3}{2} \int_1^{k+1}\mathrm{log}^2\left(\frac{k+1+y}{1+y} \right) \frac{dy}{y}.    \end{split} \]
   By Lemma \ref{htri}, it follows that
     \[
   \begin{split}
   &\;\;\;\; \mathop{\iiint}_{[1,k+1]^3}\frac{dxdydz}{(x+y)(y+z)(z+x)}\\
     &=\frac{1}{2}\mathrm{log}^3(k+1)
   -\frac{3}{2}\left(   \frac{1}{3}\mathrm{log}^3(k+1)-\zeta(2)\mathrm{log}(k+1)+\frac{7}{2}\zeta(3) +\mathrm{O}\left(\frac{\mathrm{log}\,k}{k}\right)  \right) \\
   &=\frac{3}{2}\zeta(2)\mathrm{log}\,(k+1)-\frac{21}{4}\zeta(3)+\mathrm{O}\left(\frac{\mathrm{log}\,k}{k}\right) .   \\
     \end{split} \] 
    $\hfill\Box$\\   
    
    Now we are ready to prove Theorem \ref{div}.\\
     \noindent{\bf Proof of Theorem \ref{div}:}    By Theorem \ref{sum} and Theorem \ref{tri}, one has 
     \[
     \begin{split}
     &\;\;\;\;f(k)\\
     &=\sum_{1\leq m_i\leq k}\frac{1}{(m_1+m_2)(m_2+m_3)(m_1+m_3)}-\mathop{\int}_{[1,k+1]^3}\frac{dxdydz}{(x+y)(y+z)(z+x)}\\
     &=\left[\frac{3}{2}\zeta(2) \left(\sum_{1\leq m\leq k}\frac{1}{m} \right)-\frac{41}{8}\zeta(3)+ \mathrm{O}\left( \frac{\mathrm{log}^3\;k}{k}\right)\right]\\
     &\;\;\;\;\;\;\;\;\;\;\;\;-\left[  \frac{3}{2}\zeta(2)\mathrm{log}\,(k+1)-\frac{21}{4}\zeta(3)+\mathrm{O}\left(\frac{\mathrm{log}\,k}{k}\right)            \right]\\
     &=\frac{3}{2}\zeta(2)\left( \sum_{1\leq m\leq k}\frac{1}{m} -  \mathrm{log}\,(k+1)         \right)+\frac{1}{8}\zeta(3)+ \mathrm{O}\left( \frac{\mathrm{log}^3\;k}{k}\right).
          \end{split}
     \]
     So we have $ \lim\limits_{k\to +\infty} f(k)=\frac{3}{2}\zeta(2)\gamma+\frac{1}{8}\zeta(3).$   $\hfill\Box$\\

 \section*{Acknowledgements}
 The author wants to thank Ce Xu for helpful information about special values of polylogarithms.

\end{document}